\documentclass[10pt]{amsart}

\title[Computing Twisted Gromov--Witten Invariants]{Computing Genus-Zero Twisted Gromov--Witten Invariants}

\author[Coates]{Tom Coates}
\address{Department of Mathematics\\
Imperial College London\\
180 Queen's Gate\\
London SW7 2AZ 
\\UK}
\email{t.coates@imperial.ac.uk}

\author[Corti]{Alessio Corti}
\address{Department of Mathematics\\
Imperial College London\\
180 Queen's Gate\\
London SW7 2AZ\\
UK}
\email{a.corti@imperial.ac.uk}

\author[Iritani]{Hiroshi Iritani}
\address{Faculty of Mathematics\\
Kyushu University\\
6-10-1, Hakozaki \\
Higashiku, Fukuoka, 812-8581 \\
Japan}
\curraddr{Department of Mathematics\\
Imperial College London\\
180 Queen's Gate\\
London SW7 2AZ\\
UK}
\email{iritani@math.kyoto-u.ac.jp}

\author[Tseng]{Hsian-Hua Tseng}
\address{Department of Mathematics\\
University of British Columbia\\
1984 Mathematics Road\\
Vancouver, B.C. V6T 1Z2\\
Canada}
\curraddr{Department of Mathematics\\
University of Wisconsin--Madison\\
Van Vleck Hall, 480 Lincoln Drive\\
Madison, WI 53706-1388\\
USA}
\email{tseng@math.wisc.edu}

\usepackage{amsrefs}
\usepackage{amsfonts}
\usepackage{amssymb}
\usepackage{amscd}
\usepackage{array}
\usepackage{subfigure}

\allowdisplaybreaks[3]

\newcommand{\PP}{\mathbb{P}}
\newcommand{\CC}{\mathbb{C}}
\newcommand{\ZZ}{\mathbb{Z}}
\newcommand{\RR}{\mathbb{R}}
\newcommand{\QQ}{\mathbb{Q}}
\newcommand{\Cstar}{\CC^\times}

\newcommand{\cO}{\mathcal{O}}
\newcommand{\ev}{\mathrm{ev}}
\newcommand{\bt}{\mathbf{t}}

\renewcommand{\(}{\left(}
\renewcommand{\)}{\right)}
\newcommand{\fun}{\mathbf{1}}
\DeclareMathOperator{\age}{age}

\newcommand{\cC}{\mathcal{C}}
\newcommand{\cD}{\mathcal{D}}
\newcommand{\cF}{\mathcal{F}}
\newcommand{\cH}{\mathcal{H}}
\newcommand{\cI}{\mathcal{I}}
\newcommand{\bk}{\boldsymbol{k}}

\newcommand{\cL}{\mathcal{L}}
\newcommand{\cM}{\mathcal{M}}
\newcommand{\cU}{\mathcal{U}}
\newcommand{\cX}{\mathcal{X}}
\newcommand{\cY}{\mathcal{Y}}

\newcommand{\cIX}{\mathcal{IX}}

\newcommand{\HorbX}{H^\bullet_{\text{\rm orb}}(\cX;\CC)}
\newcommand{\HorbXL}{H^\bullet_{\text{\rm orb}}(\cX;\Lambda)}
\newcommand{\HorbTX}{H^\bullet_{T,\text{\rm orb}}(\cX;\CC)}
\newcommand{\HorblocTX}{H(\cX)}

\newcommand{\HTY}{H^\bullet_T(Y;\CC)}
\newcommand{\HlocTY}{H(Y)}

\newcommand{\correlator}[1]{\left \langle #1 \right \rangle}

\DeclareMathOperator{\Eff}{Eff}
\newcommand{\EffX}{\Eff(\cX)}

\DeclareMathOperator{\Nett}{NETT}
\newcommand{\NETTX}{\Nett(\cX)}

\newcommand{\VHS}{V${\infty \over 2}$HS\ }

\DeclareMathOperator{\ch}{ch}
\newcommand{\fij}{f^{(j)}_i}

\newcommand{\Omegas}{\Omega^\tw}
\newcommand{\cLs}{\cL^\tw}

\newcommand{\bs}{{\boldsymbol{s}}}
\newcommand{\bc}{{\boldsymbol{c}}}
\newcommand{\be}{{\boldsymbol{e}}}

\newcommand{\CharV}{\mathfrak{V}}
\newcommand{\rank}{\operatorname{rank}}
\DeclareMathOperator{\Aut}{Aut}
\newcommand{\OPpair}[2]{\(#1,#2\)_{\text{\rm orb}}}
\newcommand{\CHpair}[2]{\int_{#2} #1}

\def\pair#1#2{\langle #1,#2\rangle}
\newcommand{\f}{\mathfrak{f}}

\newcommand{\fr}[1]{\left\langle #1 \right\rangle}
\newcommand{\fl}[1]{\left\lfloor #1 \right\rfloor}

\theoremstyle{plain}

\newtheorem{thm}{Theorem}[section]
\newtheorem{pro}[thm]{Proposition}
\newtheorem{lem}[thm]{Lemma}
\newtheorem*{lem*}{Lemma}

\newtheorem{cor}[thm]{Corollary}

\theoremstyle{definition}

\newtheorem{dfn}[thm]{Definition}

\newtheorem{rem}[thm]{Remark}
\newtheorem*{rem*}{Remark}
\newtheorem{exa}[thm]{Example}
\newtheorem*{exa*}{Example}
\newtheorem*{exait*}{\rm \em Example}
\newtheorem*{exadefit*}{\rm \em Example/Definition}
\newtheorem*{cla*}{\rm \em Claim}
\newtheorem*{exaA}{Example A}
\newtheorem*{exaB}{Example B}
\newtheorem*{dfn*}{Definition}

\usepackage{amsxtra} 
\newcommand{\C}{\mathbb{C}}
\newcommand{\R}{\mathbb{R}} 
\newcommand{\Z}{\mathbb{Z}} 

\newcommand{\bp}{\boldsymbol{p}} 

\newcommand{\ov}{\overline} 

\newcommand{\Spf}{\operatorname{Spf}} 
\newcommand{\Hom}{\operatorname{Hom}}

\def\parfrac#1#2{\frac{\partial #1}{\partial #2}}
\def\corr#1{\left\langle #1 \right\rangle}
\def\corrr#1{\left\langle\!\!\left\langle #1 \right\rangle\!\!\right\rangle}

\newcommand{\orb}{\text{\rm orb}}
\newcommand{\tw}{\text{\rm tw}}
\newcommand{\untw}{\text{\rm un}}
\def\<{\left\langle}
\def\>{\right\rangle}
\newcommand{\ibar}{\bar{\imath}}

\begin{document}

\maketitle

\begin{abstract}
  Twisted Gromov--Witten invariants are intersection numbers in moduli
  spaces of stable maps to a manifold or orbifold $\cX$ which depend
  in addition on a vector bundle over $\cX$ and an invertible
  multiplicative characteristic class.  Special cases are closely
  related to local Gromov--Witten invariants of the bundle, and to
  genus-zero one-point invariants of complete intersections in $\cX$.
  We develop tools for computing genus-zero twisted Gromov--Witten
  invariants of orbifolds and apply them to several examples.  We
  prove a ``quantum Lefschetz theorem'' which expresses genus-zero
  one-point Gromov--Witten invariants of a complete intersection in
  terms of those of the ambient orbifold $\cX$.  We determine the
  genus-zero Gromov--Witten potential of the type $A$ surface
  singularity $\left[\CC^2/\ZZ_n\right]$.  We also compute some
  genus-zero invariants of $\left[\CC^3/\ZZ_3\right]$, verifying
  predictions of Aganagic--Bouchard--Klemm.  In a self-contained
  Appendix, we determine the relationship between the quantum
  cohomology of the $A_n$ surface singularity and that of its crepant
  resolution, thereby proving the Crepant Resolution Conjectures of
  Ruan and Bryan--Graber in this case.
\end{abstract}

\section{Introduction}

Gromov--Witten invariants of a manifold or orbifold $\cX$ are
integrals
\begin{equation}
  \label{eq:GW}
  \int_{\left[\cX_{g,n,d}\right]^{\text{\rm vir}}} (\cdots)
\end{equation}
of appropriate cohomology classes against the virtual fundamental
class of a moduli space $\cX_{g,n,d}$ of stable maps to $\cX$.  They
give the ``virtual number'' of genus-$g$ degree-$d$ curves in $\cX$
that carry $n$ marked points constrained to lie in certain cycles
$A_1,\ldots,A_n$ in $\cX$.  The cycles $A_1,\ldots,A_n$ determine the
integrand in \eqref{eq:GW}.  It is often useful to be able to compute
similar integrals
\begin{equation}
  \label{eq:GWEuler}
  \int_{\left[\cX_{g,n,d}\right]^{\text{\rm vir}}} (\cdots) \, \be(F_{g,n,d})
\end{equation}
which involve in addition the Euler class $\be(F_{g,n,d})$ of an
\emph{obstruction bundle} $F_{g,n,d}$ over $\cX_{g,n,d}$.

\begin{exaA}
  Let $E \to \cX$ be a vector bundle which is \emph{concave}.  This
  means that
  \[
  H^0(\cC,f^\star E) = 0
  \]
  for all stable maps $f:\cC \to \cX$ of non-zero degree.  Let $d$ be
  non-zero and let $F_{g,n,d}$ be such that the fiber at the stable
  map $f:\cC \to \cX$ is
  \[
  \left. F_{g,n,d}  \right|_{f:\cC \to \cX} = H^1(\cC,f^\star E).
  \]
  Then integrals \eqref{eq:GWEuler} are Gromov--Witten invariants of
  the (non-compact) total space of $E$: they are \emph{local
    Gromov--Witten invariants} \cite{Chiang--Klemm--Yau--Zaslow}.
\end{exaA}

\begin{exaB}
  Let $E \to \cX$ be a vector bundle which is \emph{convex}.  This
  means that 
  \[
  H^1(\cC,f^\star E) = 0
  \]
  for all genus-zero one-pointed stable maps\footnote{One-pointed
    stable maps are those with $n=1$.} $f:\cC \to \cX$.  Let
  $F_{0,1,d}$ be such that
 \[
 \left. F_{0,1,d} \right|_{f:\cC \to \cX} = H^0(\cC,f^\star E).
  \]
  Then integrals \eqref{eq:GWEuler} with $g=0$ and $n=1$ give
  Gromov--Witten invariants of a suborbifold of $\cX$ cut out by a
  section of $E$.
\end{exaB}

In this paper we consider \emph{twisted Gromov--Witten invariants}.
These are integrals 
\begin{equation}
  \label{eq:GWtwisted}
  \int_{\left[\cX_{g,n,d}\right]^{\text{\rm vir}}} (\cdots) \, \bc(F_{g,n,d})
\end{equation}
involving an invertible multiplicative characteristic class\footnote{A
  characteristic class $\bc$ is multiplicative if $\bc(E_1 \oplus E_2)
  = \bc(E_1) \bc(E_2)$.  It is invertible if $\bc(E)$ is invertible in
  $H^\bullet(\cY)$ whenever $E$ is a vector bundle over $\cY$.
  Invertible multiplicative characteristic classes extend to
  $K$-theory: $\bc(E_1 \ominus E_2) = \bc(E_1) \bc(E_2)^{-1}$.}  $\bc$
applied to an ``obstruction $K$-class'' $F_{g,n,d} \in
K^0\(\cX_{g,n,d}\)$,
\[
\left. F_{g,n,d} \right|_{f:\cC \to \cX} = H^0(\cC,f^\star F) \ominus
H^1(\cC,f^\star F),
\]
where $F$ is a vector bundle over $\cX$.  (We give a formal definition
in Section~\ref{sec:twisted} below.)  When $\bc$ is the trivial
characteristic class, these coincide with ordinary Gromov--Witten
invariants.  The Euler class is not invertible, but nonetheless
Examples A and B can be included in this framework as follows.  Every
vector bundle $F$ carries the action of a torus $T$ which rotates
fibers and leaves the base invariant; we can always take $T = \Cstar$,
and if $F$ is the direct sum of line bundles then we can take $T =
\(\Cstar\)^{\rank F}$.  The $T$-equivariant Euler class is invertible
over the fraction field of $H_T^\bullet\(\{\text{pt}\}\)$.  Example B
arises by taking $F=E$ and $\bc$ to be the $T$-equivariant Euler
class, and then taking the non-equivariant limit.  Example A arises by
taking $F=E$ and $\bc$ to be the $T$-equivariant inverse Euler class,
and then taking the non-equivariant limit.  Twisted Gromov--Witten
invariants also occur in virtual localization formulas
\cite{Graber--Pandharipande} for the $T$-equivariant Gromov--Witten
invariants of an orbifold $\cY$ equipped with the action of a torus
$T$.  There $\cX_{g,n,d}$ is part of the $T$-fixed substack of the
moduli stack of stable maps to $\cY$ and $\bc$ is the $T$-equivariant
inverse Euler class.  If we can compute twisted Gromov--Witten
invariants, therefore, then we can compute local Gromov--Witten
invariants, genus-zero one-point invariants of complete intersections,
and $T$-equivariant Gromov--Witten invariants.  Twisted Gromov--Witten
invariants for other choices of $\bc$ can be interpreted as
Gromov--Witten invariants with values in generalized cohomology
theories \cite{Givental:symplectic}.

When $\cX$ is a manifold, one can compute twisted Gromov--Witten
invariants using results of Coates--Givental
\cite{Coates--Givental:QRRLS}.  They prove a ``quantum Riemann--Roch
theorem'' expressing twisted Gromov--Witten invariants of all genera,
for any choice of $\bc$ and $F$, in terms of ordinary Gromov--Witten
invariants of $\cX$.  From this they deduce a ``quantum Lefschetz
theorem'' which gives simple closed formulas for genus-zero twisted
invariants in the case where $F$ is the direct sum of convex line
bundles and $\bc$ is the $T$-equivariant Euler class.  This implies
most of the known mirror theorems for toric complete intersections.
The results in \cite{Coates--Givental:QRRLS} are based on a
Grothendieck--Riemann--Roch argument, essentially due to Mumford
\cite{Mumford} and Faber--Pandharipande \cite{Faber--Pandharipande},
and a geometric formalism introduced by Givental
\cite{Givental:quantization}.  A quantum Riemann--Roch theorem for
orbifolds has been established by Tseng \cite{Tseng} using the
Grothendieck--Riemann--Roch theorem of Toen \cite{Toen}. Tseng also
proved a version of quantum Lefschetz in the orbifold setting
\cite[Theorem 5.15]{Tseng}, but this holds only under very restrictive
hypotheses on the bundle $F$.

In this paper we prove a much more general quantum Lefschetz-style
result for orbifolds.  This is Theorem~\ref{thm:generaltwist} below.
It applies whenever $F$ is a direct sum of line bundles, without
restriction on the invertible multiplicative class $\bc$, and
determines genus-zero twisted Gromov--Witten invariants of an orbifold
$\cX$ in terms of the ordinary Gromov--Witten invariants of $\cX$.  It
removes many of the restrictive hypotheses from Tseng's result and (as
did \cite{Iritani--V3}) improves on Coates--Givental when $\cX$ is a
manifold, in that:
\begin{itemize}
\item the characteristic class $\bc$ does not have to be an Euler
  class; and
\item the bundle $F$ is not assumed to be convex.
\end{itemize}

In practice Theorem~\ref{thm:generaltwist} is most useful in the
situation of Examples A and B.  This gives nothing new in the manifold
setting, as these cases were already covered by
\cite{Coates--Givental:QRRLS}, but the improvement for orbifolds is
significant.  We illustrate this with several examples.  In
Section~\ref{sec:hypersurfaces} we consider the situation of
Example~A, computing certain genus-zero Gromov--Witten invariants of
the quintic hypersurface in $\PP(1,1,1,1,2)$.  We also prove a quantum
Lefschetz theorem for orbifolds, Corollary~\ref{thm:orbQL} below,
which directly generalizes \cite[Theorem 2]{Coates--Givental:QRRLS}
and \cite[Theorem 5.15]{Tseng}.  (This suffices, for example, to
determine the even-degree part of the small quantum orbifold
cohomology algebra of any of the $181$ Fano $3$-fold weighted
projective complete intersections with terminal singularities: see
\cite[Proposition~1.10]{CCLT}.)  In Section~\ref{sec:local} we
consider the situation of Example~B, computing in Section~\ref{sec:An}
the genus-zero Gromov--Witten potential of the type $A$ surface
singularity $\left[\CC^2/\ZZ_n\right]$.  This has been determined for
$n=2$ by Bryan--Graber \cite{Bryan--Graber}; for $n=3$ by
Bryan--Graber--Pandharipande \cite{Bryan--Graber--Pandharipande}; for
$n=4$ by Bryan--Jiang \cite{Bryan--Jiang}.  Their methods are quite
different from ours.  Perroni \cite{Perroni} has studied the small
quantum cohomology of orbifolds with transverse $ADE$ singularities,
and part of the potential for $\left[\CC^2/\ZZ_n\right]$ can be
extracted from his results.  Maulik \cite{Maulik} has computed the
genus-zero Gromov-Witten potential and certain higher-genus
Gromov-Witten invariants of $[\CC^2/\ZZ_n]$.  In
Section~\ref{sec:C3Z3} we compute certain genus-zero Gromov--Witten
invariants of $\left[\CC^3/\ZZ_3\right]$ where $\ZZ_3$ acts with
weights $(1,1,1)$, verifying predictions of Aganagic--Bouchard--Klemm
\cite{ABK}. In Appendix~\ref{sec:CRC}, which can be read separately
from the main text, we combine results from Section~\ref{sec:An} with
arguments from toric mirror symmetry to prove the Crepant Resolution
Conjectures of Ruan and Bryan--Graber for the type $A$ surface
singularity $\left[\CC^2/\ZZ_n\right]$: this is new for $n \geq 5$.
In Appendix~\ref{sec:rigorouscone} we prove some foundational results,
describing certain aspects of Givental's geometric formalism in terms
of non-Noetherian formal schemes.

\subsection*{Acknowledgments} We are grateful to Yongbin Ruan for many
productive and inspiring conversations, and to Jim Bryan, Etienne
Mann, and Fabio Perroni for useful discussions.  We are grateful also
to the anonymous referees for their careful reading and valuable
comments.  In particular, Appendix~\ref{sec:rigorouscone} was added at
the suggestion of one of the referees.  T.C. thanks Bong Lian and
Shing-Tung Yau for helpful conversations.  H.I. is grateful to Akira
Ishii for teaching him about autoequivalences of the derived category.
This work was begun while T.C., H.I., and H.-H.T. held postdoctoral
fellowships at the Mathematical Sciences Research Institute as part of
the ``New Topological Structures in Physics'' program.  T.C. was in
addition supported by the Royal Society and by NSF grant DMS-0401275.
H.I. was in addition supported by the Grant-in-Aid for Scientific
Research 18-15108 and the 21st Century COE program of Kyushu
University.  H.-H.T. was in addition supported by Institut
Mittag-Leffler.

\section{Preliminaries}

In this section we fix notation for orbifold cohomology and orbifold
Gromov--Witten theory.  These notions were introduced by Chen--Ruan in
the symplectic category; an algebraic version of the theory has
been developed by Abramovich--Graber--Vistoli.  We will assume that the
reader is familiar with this material --- see \cite[Section 2]{CCLT}
for a brief overview and the original sources
\cite{Chen--Ruan:orbifold,Chen--Ruan:GW, AGV:2} for a comprehensive
treatment.

\subsection{Orbifold Cohomology}

We work in the algebraic category, using notation as follows.

\begin{center}
  \begin{tabular}{lp{0.75\textwidth}} $\cX$ & a proper smooth
    Deligne--Mumford stack over $\CC$ with projective coarse moduli space. \\
    $\cIX$ & the inertia stack of $\cX$.  A point of $\cIX$ is a pair
    $(x,g)$ with $x$ a point of $\cX$ and $g \in \Aut_{\cX}(x)$.\\
    $\cIX = \coprod_{i \in \cI} \cX_i$ & the decomposition of $\cIX$
    into components; here $\cI$ is an index set \\
    $I$ & the involution of $\cIX$ which sends $(x,g)$ to
    $(x,g^{-1})$. \\
    $\HorbX$ & the orbifold cohomology groups of $\cX$.  These are the
    cohomology groups $H^\bullet(\cIX;\CC)$ of the inertia stack.  \\
    $\age$ & a rational number associated to each component $\cX_i$ of
    the inertia stack. Chen--Ruan call this the \emph{degree-shifting number}. \\
    $\OPpair{\alpha}{\beta}$ & the
    orbifold Poincar\'e pairing $\int_{\cIX} \alpha \cup I^\star \beta$ \\
  \end{tabular}
\end{center}

The grading on orbifold cohomology is shifted by the age: $\alpha\in
H^p(\cX_i;\CC)$ has degree $\deg\alpha=p + 2 \age(\cX_i)$.

\subsection{Moduli Spaces of Stable Maps}

Let $\cX_{0,n,d}$ denote, as in \cite[Section 2.2.1]{CCLT}, the moduli
stack of $n$-pointed genus-zero stable maps to $\cX$ of degree $d \in
H_2(\cX;\QQ)$.  This is almost exactly what
Abramovich--Graber--Vistoli call the \emph{stack of twisted stable
  maps} $\mathcal{K}_{0,n}(\cX,d)$.  The only difference is that they
regard the degree as a curve class on the coarse moduli space of
$\cX$, whereas we regard it as an element of $H_2(\cX;\QQ)$.  We will
not use the term ``twisted stable map'' as for us ``twisted'' means
something different.

There are \emph{evaluation maps} $\ev_i:\cX_{0,n,d}\to
\overline{\cI\cX}$, one for each marked point, which take values in the
rigidified cyclotomic inertia stack $\overline{\cIX}$.  Since there
is a proper \'etale surjection $\cIX \to \overline{\cIX}$, we can use
the evaluation maps to define cohomological pull-backs
\begin{align*}
(\ev_i)^\star &: \HorbX \to H^\bullet(\cX_{0,n,d};\CC) 
\end{align*}
even though the maps $\ev_i$ do not take values in the inertia stack
$\cIX$.  We write
\[
[\cX_{0,n,d}]^{\text{\rm vir}}\in H_\bullet(\cX_{0,n,d};\CC)
\]
for the virtual fundamental class of the moduli stack and 
\begin{align*}
\psi_i \in H^2(\cX_{0,n,d};\CC), && i \in \{1,2,\ldots,n\},
\end{align*}
for the first Chern class of the universal cotangent line bundle
$L_i$. The fiber of $L_i$ at the stable map $f : \cC\to \cX$ is the
cotangent line to the {\em coarse moduli space of $\cC$} at the $i$th
marked point. 

\subsection{Twisted and Untwisted Gromov--Witten Invariants} 
\label{sec:twisted}

A more detailed account of the material in this section can be found
in \cite{Tseng}.  Twisted Gromov--Witten invariants are a
family of invariants of $\cX$ which depend on an invertible
multiplicative characteristic class $\bc$ and a vector bundle $F \to
\cX$.  Throughout this paper we will take $F$ to be the direct sum of
line bundles,
\[
F = \bigoplus_{j=1}^{j=r} F^{(j)}.
\]
In applications below we will need to take $\bc$ to be a
$T$-equivariant cohomology class, where the torus $T=(\Cstar)^{r}$ acts
on $F$ by scaling the fibers.  We write
\[
H_T^\bullet(\{{\rm pt}\}) = \CC[\lambda_1,\ldots,\lambda_r]
\]
where $\lambda_i$ is Poincar\'e-dual to a hyperplane in the $i$th
factor of $\(\CC \mathbb{P}^\infty\)^r \cong BT$.

Consider the universal family over $\cX_{0,n,d}$
 \[
\begin{CD}
   \cC_{0,n,d} @> f>> \cX \\ @V{\pi}VV \\ \cX_{0,n,d}
\end{CD} 
\]
and define an element $F_{0,n,d} \in K^0(\cX_{0,n,d})$ by 
\[
F_{0,n,d}:=\pi_! f^\star F,
\]
where $\pi_!$ is the $K$-theoretic push-forward. Genus-zero
\emph{twisted Gromov--Witten invariants of $\cX$} are intersection
numbers of the form
\begin{equation}
  \label{eq:correlator}
  \<\alpha_1\psi^{k_1},\ldots,\alpha_n \psi^{k_n}\>_{0,n,d}^{\cX,\tw} := 
  \int_{[\cX_{0,n,d}]^{\text{\rm{vir}}}}\bc(F_{0,n,d})\cup\prod_{i=1}^n \ev_i^\star(\alpha_i)
  \cdot \psi_i^{k_i}
\end{equation}
where $\alpha_1,\ldots,\alpha_n \in \HorbX$; $k_1,\ldots,k_n$ are
non-negative integers; and the integral denotes cap product with the
virtual fundamental class.  If $\bc$ is the trivial characteristic
class --- this is the case of usual, untwisted Gromov--Witten
invariants --- then we will replace the superscript ``$\tw$'' by
``$\untw$''.

\begin{rem}\label{remarks_on_twisted_inv} The Gromov--Witten
  invariants defined here coincide with those considered in
  \cite{Tseng}: we use slightly different stacks of stable maps and
  also a different definition of the pull-back $(\ev_i)^\star$, but
  these two differences cancel each other out.  The descendant class
  denoted $\psi_i$ here is denoted in \cite{Tseng} by $\bar{\psi_i}$.
\end{rem}
 
Genus-zero twisted orbifold Gromov--Witten invariants together define
a Frobenius manifold, as we now explain.  Fix a K\"ahler class
$\omega$ on $\cX$.  Let $\EffX$ be the semigroup of degrees of
representable maps from possibly-stacky curves to $\cX$ (\emph{i.e.}
of degrees of effective curves in $\cX$) and define the \emph{Novikov
  ring} $\Lambda$ to be the completion of the group ring $\CC[\EffX]$
of $\EffX$ with respect to the additive valuation $v$,
\[
v\(\sum_{d\in \EffX} a_d Q^d\) = \min_{a_d\neq 0}\CHpair{\omega}{d}, 
\]
where $Q^d$ is the element of $\CC[\EffX]$ corresponding to $d\in
\EffX$.  Note that the completion depends on the choice of K\"{a}hler
class $\omega$.  The Frobenius manifold is based on the free
$\Lambda$-module
\[
\HorbXL :=H^\bullet(\cIX;\CC)\otimes_\CC \Lambda.
\] 
To define the pairing, observe that the inertia stack $\cI F$ of the
total space of the vector bundle $F \to \cX$ is a vector bundle over
$\cIX$ --- the fiber of $\cI F$ over the point $(x,g) \in \cIX$
consists of the $g$-fixed subspace of the fiber of $F$ over $x$ ---
and set
\begin{equation}
  \label{eq:twistedpairing}
  \(\alpha,\beta\)^{\tw}_\orb = \int_{\cIX} \alpha \cup
  I^\star \beta\cup \bc\(\cI F\).
\end{equation}

\begin{exa}
  Let $\cX = B \mu_r$ and let $F \to \cX$ be the tautological line
  bundle.  Then $\cIX$ is the disjoint union of $r$ copies of $B\mu_r$
  where the $j$th copy, $0 \leq j < r$, corresponds to the element
  $\zeta^j = \exp\({2 \pi \sqrt{-1} j \over r}\) \in \mu_r$.  For $j
  \ne 0$, $\zeta^j$ acts non-trivially on the fiber of $F$ and so the
  fiber of $\cI F$ over the $j$th copy of $B\mu_r$ is the
  zero-dimensional vector space.  The restriction of $\cI F$ to the
  zeroth copy of $B\mu_r$ is $F$.
\end{exa}

Genus-zero twisted Gromov--Witten invariants assemble to give a family
of products, defined by
\begin{equation}
  \label{eq:twistedmultiplication}
  \(\alpha\bullet_\tau \beta,\gamma\)^{\tw}_\orb =
  \sum_{d\in\EffX}\sum_{n\geq
    0}\frac{Q^d}{n!}\<\alpha,\beta,\gamma, \tau, \tau,\ldots,\tau
  \>_{0,n+3,d}^{\cX,\tw},
\end{equation}
parametrized by $\tau$ in a formal neighbourhood of zero in $\HorbXL$.
When $\bc=1$, this gives the usual Frobenius manifold structure on
orbifold cohomology.

\section{Givental's Symplectic Formalism} 
\label{sec:GiventalFormalism} 

In this section we will descibe how to encode genus-zero twisted
orbifold Gromov--Witten invariants in a Lagrangian submanifold of a
certain symplectic vector space.  This idea is due to Givental
\cite{Givental:symplectic}; it was adapted to the orbifold setting by
Tseng \cite{Tseng}.  We will describe only the aspects of the theory
which we need, referring the reader to \cite{Givental:symplectic,
  Tseng} and the references therein for motivation, context, and
further examples of this approach.  In particular the genus-zero
picture used here is only part of a more powerful formalism involving
Gromov--Witten invariants of all genera, and we will not discuss this
at all.

\begin{dfn*}
  For a topological ring $R$ with a non-negative 
  additive valuation $v:R\setminus \{0\} 
  \rightarrow \RR_{\ge 0}$, 
  define the space of \emph{convergent Laurent series} in $z$ to be  
  \begin{align*}
    R\{z,z^{-1}\} & := \left\{\sum_{n\in \ZZ} r_n z^n\; :\; r_n\in R,
      \ v(r_n) \to \infty \text{ as } |n|\to \infty \right\}.
    \intertext{If $R$ is complete, this becomes a
      ring\footnotemark. Set} R\{z\} & := \left\{\sum_{n\geq 0} r_n
      z^n\; :\; r_n\in R, \ v(r_n) \to \infty \text{ as } n\to \infty
    \right\}.
  \end{align*}
  \footnotetext{In this case, $R\{z,z^{-1}\}$ coincides with the
    completion of $R[z,z^{-1}]$ under the induced valuation with
    $v(z)=0$.}
\end{dfn*}

Consider the space of orbifold-cohomology-valued convergent Laurent
series
\[
\cH:=\HorbX \otimes \Lambda\{z,z^{-1}\} 
\]
equipped with the $\Lambda$-valued symplectic form
\[
\Omegas(f,g):=\text{Res}_{z=0}\(f(-z),g(z)\)^\tw_\orb \, dz.
\]
We encode genus-zero twisted orbifold Gromov--Witten invariants via
the germ of a Lagrangian submanifold $\cLs$ of $(\cH,\Omegas)$,
defined as follows.  Let 
\begin{align*}
  \{\phi_\alpha : 1 \leq \alpha \leq N\} && \text{and} && \{\phi^\alpha : 1 \leq \alpha
\leq N\}
\end{align*}
be $\Lambda$-bases for $\HorbXL$ which are dual with repect to the
pairing \eqref{eq:twistedpairing}.  The submanifold-germ $\cLs$
consists of all points of $\cH$ of the form
\begin{multline}
  \label{eq:pointonL}
  -z + t_0 + t_1 z + t_2 z^2 + \cdots \\
  + \sum_{\substack{d \in \EffX \\ n \geq 0}} 
  \sum_{\substack{ i_1,\ldots,i_n\\\alpha_1,\ldots,\alpha_n}} 
  \sum_{\substack{k \geq 0 \\ 1 \leq \epsilon \leq N}}
  {Q^d t^{\alpha_1}_{i_1} \cdots t^{\alpha_n}_{i_n} \over n!}
  \< \phi_{\alpha_1}
  \psi^{i_1},\ldots,\phi_{\alpha_n} \psi^{i_n},\phi_\epsilon \psi^k
  \>^{\cX,\tw}_{0,n+1,d}
  { \phi^\epsilon \over (-z)^{k+1}}
\end{multline}
where $t_0 + t_1 z + t_2 z^2 + \ldots$ lies in a formal neighbourhood of
zero in $\HorbX \otimes \Lambda\{z\}$ and $t_i = \sum_{\alpha} t_i^\alpha
\phi_\alpha$.  If we write $\bt(z) = t_0 + t_1 z + t_2 z^2 + \ldots$
then \eqref{eq:pointonL} is
\[
  -z + \bt(z) 
  + \sum_{n \geq 0} \sum_{d \in \EffX} 
  \sum_{1 \leq \epsilon \leq N}
  {Q^d \over n!}
  \< \bt(\psi),\ldots,\bt(\psi),{\phi_\epsilon \over -z-\psi}  \>^{\cX,\tw}_{0,n+1,d}
  \phi^\epsilon.
\]
The submanifold-germ $\cLs$ has extremely special geometric
properties, which are listed in 
\cite[Theorem 2.15]{CCIT:wallcrossings1} and
\cite[Section 3.1]{Tseng}. 
These follow from the fact that genus-zero
twisted Gromov--Witten invariants satisfy the String Equation, the
Dilaton Equation, and the Topological Recursion Relations
\cite{Givental:symplectic, Tseng}. 

\subsection*{A Remark on Rigour} The definition of $\cLs$ just given
is not completely rigorous, as we did not spell out exactly what we
mean by a formal neighbourhood in an infinite-dimensional vector
space.  A rigorous definition of $\cLs$, as a non-Noetherian
\emph{formal scheme}, is given in Appendix \ref{sec:rigorouscone}.
There we also establish various geometric properties of $\cLs$ which
will be needed later: see Propositions~\ref{prop:cone},
\ref{prop:closedunderz}, \ref{prop:tangentspace_derivativesofJ}, and
Corollary~\ref{cor:Dmod}.  The rest of this paper can therefore be
read in two ways.  The reader who is happy to work with an intuitive
notion of formal neighbourhood can read the rest of the text as it is,
omitting Appendix~\ref{sec:rigorouscone}.  The discussion will
then be informal, but no serious confusion should result.  The reader
who prefers a completely formal approach should at this point skip to
Appendix~\ref{sec:rigorouscone}, and replace the definition
\eqref{eq:pointonL} above with definition \eqref{eq:def_rigorouscone}
below.  The rest of the text can then be read as a series of rigorous
arguments within the framework constructed in
Appendix~\ref{sec:rigorouscone}.

\subsection{The Twisted $J$-Function}
\label{sec:twistedJ}

Let 
\begin{equation}
  \label{eq:twistedJfun}
  J^\tw(\tau,z) = z + \tau + 
  \sum_{n \geq 0} \sum_{d \in \EffX} 
  \sum_{1 \leq \epsilon \leq N}
  {Q^d \over n!}
  \< \tau,\tau,\ldots,\tau,{\phi_\epsilon \over z-\psi}  \>^{\cX,\tw}_{0,n+1,d}
  \phi^\epsilon.
\end{equation}
This formal power series in the components $\tau^1,\ldots,\tau^N$ of
$\tau = \tau^1 \phi_1 + \ldots + \tau^N \phi_N$, called the
\emph{twisted $J$-function of $\cX$}, will play an important role
below.  It takes values in $\cH$ and gives a distinguished family
\begin{align*}
  \tau \longmapsto J^\tw(\tau,-z), && \text{$\tau$ in a formal
    neighbourhood of zero in $\HorbXL$,}
\end{align*}
of elements of $\cLs$ characterized among such families by the
property that
\begin{equation}
  \label{eq:Jcharact}
  J^\tw(\tau,-z) = -z + \tau + O(z^{-1}).
\end{equation}

We write $J^\untw$, $\cL^\untw$, and $\Omega^\untw$ for the
specializations of, respectively, $J^\tw$, $\cL^\tw$, and $\Omega^\tw$
to the case $\bc = 1$ --- \emph{i.e.} for the corresponding objects in
untwisted Gromov--Witten theory.  The untwisted $J$-function satisfies
a system of differential equations
\begin{equation}
  \label{eq:qdes}
  z {\partial \over \partial \tau^\alpha} {\partial \over \partial
    \tau^\beta} J^\untw(\tau,z) = \sum_{\gamma=1}^N c_{\alpha
    \beta}^{\phantom{\alpha \beta} \gamma}(\tau) {\partial
    \over \partial \tau^\gamma} J^\untw(\tau,z)
\end{equation}
where $c_{\alpha \beta}^{\phantom{\alpha \beta} \gamma}(\tau)$ are the
structure constants of the untwisted multiplication with respect to
the basis $\{\phi_\epsilon\}$:
\[
\phi_\alpha \bullet_\tau \phi_\beta \, \bigg|_{\bc = 1} =
\sum_{\gamma=1}^N c_{\alpha \beta}^{\phantom{\alpha \beta}
  \gamma}(\tau) \phi_\gamma.
\]
One can see this either as a consequence of the geometric properties
of $\cL^\untw$ \cite{Givental:symplectic}, or directly from the
Topological Recursion Relations as in \cite[Lemma~2.4]{CCLT} or
\cite[Proposition~2]{Pandharipande}.

\section{General Twists in Genus Zero}

In this section we give a formula for a family of elements on the
Lagrangian submanifold $\cLs$ for the twisted theory.  The key
ingredient is Tseng's genus-zero orbifold quantum Riemann--Roch
theorem \cite{Tseng}, so we begin by stating this.

\subsection{Orbifold Quantum Riemann--Roch}
\label{sec:QRR}
Given a line bundle $L \to \cX$ and a geometric point $(x,g) \in
\cIX$, there is a unique rational number $f \in [0,1)$ such that $g$
acts on the fiber of $L$ over $x$ by multiplication by $\exp\(2 \pi
\sqrt{-1} f\)$.  The value of $f$ depends only on the component
$\cX_i$ of $\cIX$ containing $(x,g)$.  Since $F$ is the direct sum of
line bundles,
\[
F = \bigoplus_{j=1}^{j=r} F^{(j)},
\]
this defines a collection of rational numbers $f_i^{(j)}$ where $1
\leq j \leq r$ and $i \in \cI$.  The other ingredients in the
statement are the first Chern classes $\rho^{(j)} \in H^2(\cX;\CC)$ of
$F^{(j)}$, regarded as elements of orbifold cohomology via the natural
inclusion $\cX \to \cIX$, and the (unique) sequence of parameters
$s_0,s_1,s_2,\ldots$ such that
\begin{equation}
\label{eq:universalclass}
\bc(\cdot) = \exp\(\sum_{k \geq 0} s_k \ch_k(\cdot) \).
\end{equation} 
Here $\ch_k$ is the $k$th component of the Chern character.  We add
the variables $s_k$ to our ground ring, working henceforth over the
completion $\Lambda[\![s_0,s_1,\ldots]\!] $ of
$\CC[\EffX][s_0,s_1,\dots]$ with respect to the additive valuation $v$
such that
\[
v(Q^d)=\CHpair{\omega}{d}, \quad v(s_k)=k+1.   
\]
Later we will need the notation
\[
\bs(x) = \sum_{k \geq 0} s_k {x^k \over k!}.
\]
Recall that the Bernoulli polynomials $B_n(x)$ are defined by
\[ 
\sum_{n=0}^\infty B_n(x)\frac{z^n}{n!} = \frac{ze^{zx}}{e^z-1}.
\] 

\begin{thm}[{\cite[Corollary 1]{Tseng}}] \label{thm:genuszeroQRR} The
  transformation $\Delta^\tw\colon \cH \rightarrow \cH$ defined by
  \[ 
  \Delta^{\tw}=\bigoplus_{i\in \cI} \prod_{j=1}^r \exp\left (\sum_{l,m\ge 0}
    s_{l+m-1} \frac{B_m(\fij)}{m!} \frac{(\rho^{(j)})^{l}}{l!}
    z^{m-1}\right),
  \] 
  where $\rho^{(j)}$ acts on $\cH$ via the Chen--Ruan orbifold cup
  product \cite{Chen--Ruan:orbifold,AGV:2}, gives a linear
  symplectomorphism between $(\cH,\Omega^\untw)$ and
  $(\cH,\Omega^\tw)$, and
  \[ 
  \cL^\tw = \Delta^\tw \(\cL^\untw\).
  \]
\end{thm}
Here we are implicitly using the facts that
\[
\HorbX = \bigoplus_{i \in \cI} H^\bullet(\cX_i;\CC)
\]
and that the action of $\rho^{(j)}$ preserves this decomposition. 
The Chen-Ruan orbifold cup product by $\rho^{(j)}$ 
coincides with the ordinary cup product by $\pi^*\rho^{(j)}$ 
\cite[Lemma 2.3.7]{Tseng}, 
where $\pi\colon \cIX \to \cX$ is the natural projection.  
We define $s_{-1}$ to be zero.

\begin{rem} Multiplication by $\sqrt{\bc\(\cI F\)}$, using the usual
  cup product on $H^\bullet(\cIX;\CC)$, gives an isomorphism between
  the symplectic vector spaces $(\cH,\Omega^\tw)$ and
  $(\cH,\Omega^\untw)$. The transformation $\Delta^\tw$ appearing
  above differs from that in \cite{Tseng} because the transformation
  there was regarded as an automorphism of $(\cH,\Omega^\untw)$ via
  this identification $(\cH,\Omega^\tw) \cong (\cH,\Omega^\untw)$
\end{rem}

\subsection{A Family of Elements of $\cL^\tw$}  
\label{sec:Ifunction}

It will be convenient to break up the untwisted $J$-function
$J^\untw(\tau,z)$ into contributions from stable maps of different
topological types.

\begin{dfn}
  The \emph{topological type} of a stable map $f:\cC \to \cX$, where
  $\cC$ has genus $g$ and marked points $x_1,\ldots,x_n$ and $f$ has
  degree $d \in H_2(\cX;\QQ)$, is the triple
  \[
  \theta = (g,d,S)
  \]
  where $S$ is the ordered $n$-tuple consisting the elements of
  $\cI$ which label the components of the inertia stack picked out by
  the marked points $x_1,\ldots,x_n$.
\end{dfn}

The topological type is constant on each component of the moduli space
$\cX_{0,n,d}$.  We write $\NETTX$ for the set of all topological types
of stable maps to $\cX$, or in other words for the set of
\emph{effective topological types} in $\cX$, and $J_\theta(\tau,z)$ for
the contribution to the untwisted $J$-function from stable maps of
topological type $\theta$, so that
\[
J^\untw(\tau,z) = \sum_{\theta \in \NETTX} J_\theta(\tau,z).
\]

\begin{rem}
  \label{rem:decomposition}
  In practice one can determine this decomposition by choosing the
  basis $\{\phi_\alpha\}$ for $\HorbXL$ so that each $\phi_\alpha$ is
  supported on exactly one component $\cX_{i(\alpha)}$ of $\cIX$.
  Then the term 
  \[
  {Q^d \tau^{\alpha_1} \cdots \tau^{\alpha_n} \over n! \,
    (z)^{k+1}}  
  \< \phi_{\alpha_1},\ldots,\phi_{\alpha_n},\phi_\epsilon \psi^k
  \>^{\cX,\untw}_{0,n+1,d} \phi^\epsilon 
  \]
  in the power series expansion \eqref{eq:twistedJfun} of
  $J^\untw(\tau^1\phi_1+\cdots+\tau^N\phi_N,z)$ contributes to
  $J_\theta$ only for $\theta = (0,d,S)$ where $S = \left(
    i(\alpha_1),i(\alpha_2),\ldots,i(\alpha_n),i(\epsilon) \right)$.
 \end{rem}

 \begin{lem} \label{lem:propertiesofJ}
   Let $\theta \in \NETTX$ be the topological type $(0,d,S)$ where $S
   = (i_1,\ldots,i_n)$.  Then
   \begin{enumerate}
   \item The orbifold cohomology class $J_\theta(\tau,z)$ is supported
     on the component $I(\cX_{i_n})$ of $\cIX$.
   \item If $D_i$ is the dilation vector field on
     $H^\bullet(\cX_i;\CC) \subset \HorbX$, so
     \[ 
     D_i = \sum_{\nu} x^\nu {\partial \over \partial x^\nu}
     \] 
     for any linear co-ordinate system $(x^\nu)$ on
     $H^\bullet(\cX_i;\CC)$, then
     \[
     D_i J_\theta(\tau,z) = n_i J_\theta(\tau,z)
     \]
     where $n_i$ is the number of times that $i$ occurs in
     $(i_1,\ldots,i_{n-1})$.
   \item If $\rho \in H^2(\cX;\CC)$ is regarded as an orbifold
     cohomology class via the natural inclusion $\cX \to \cIX$ then
     \[
     z \nabla_\rho J_\theta(\tau,z) = \( \rho + \<\rho,d\>z\)
     J_\theta(\tau,z).
     \]
     Here $\nabla_\rho$ is the directional derivative: 
     $\nabla_\rho J_\theta(\tau) := \frac{d}{dt} J_\theta(\tau + t \rho)|_{t=0}$. 
   \end{enumerate}
 \end{lem}

 \begin{proof}
   (1) and (2) follow immediately from Remark~\ref{rem:decomposition}.
   (3) follows from the Divisor Equation \cite[Theorem 8.3.1]{AGV:2}.
 \end{proof}

 Consider a topological type $\theta \in \NETTX$ with $\theta =
 (0,d,S)$ and $S = (i_1,\ldots,i_n)$.  Let $\ibar_n \in \cI$ be such
 that the component $\cX_{\ibar_n}$ is $I(\cX_{i_n})$, and let
\[
N^{(j)}_\theta = \<\rho^{(j)},d\> - f_{i_1}^{(j)} - f_{i_2}^{(j)} - \cdots -
f_{i_{n-1}}^{(j)} + f_{\ibar_n}^{(j)}.
\]
Riemann--Roch for orbifold curves implies that $N^{(j)}_\theta$ is an
integer: for any stable map $h:\cC \to \cX$ of topological type
$\theta$, the Euler characteristic
\begin{align*}
  \chi\(\cC,h^\star F^{(j)}\) &=  1+\pair{\rho^{(j)}}{d} - f^{(j)}_{i_1} - f^{(j)}_{i_2}
  - \cdots - f^{(j)}_{i_n},
\intertext{and}
f_{\ibar_n}^{(j)} &=
\begin{cases}
  0 & \text{if $f_{i_n}^{(j)} = 0$} \\
  1 - f_{i_n}^{(j)} & \text{if $f_{i_n}^{(j)} \ne 0$}.
\end{cases}
\end{align*}
Define the \emph{modification factor}
\begin{equation*}
  M_\theta(z) = 
  \prod_{j=1}^{j=r}
  {
    \displaystyle \prod_{- \infty < m \leq N^{(j)}_\theta} 
    \exp\left[\bs\(\rho^{(j)}+\(m-f_{\ibar_n}^{(j)}\)z\)\right]
    \over
    \displaystyle \prod_{- \infty < m \leq 0} 
    \exp\left[\bs\(\rho^{(j)}+\(m-f_{\ibar_n}^{(j)}\)z\)\right]
  }
\end{equation*}
and set
\begin{equation*}
  I^\tw(\tau,z)  = \sum_{\theta \in \NETTX} M_\theta(z) \cdot J_\theta(\tau,z).
\end{equation*}
where the multiplication is with respect to the Chen--Ruan orbifold
cup product.  $I^\tw$ is a formal power series in the components
$\tau^1, \ldots, \tau^N$ of $\tau$ which takes values in $\cH$.  When
$\cX$ is a variety and $\bc$ is the $T$-equivariant Euler class, it
coincides with the \emph{hypergeometric modification} $I_F$ in
\cite[Section 7]{Coates--Givental:QRRLS}.

\begin{thm}
\label{thm:generaltwist} The family
\[ 
\tau \mapsto I^\tw(\tau,-z)
\] 
of elements of $\cH$ lies on the Lagrangian submanifold $\cLs$.
\end{thm}

\subsection*{Remark} In the formal framework developed in
Appendix~\ref{sec:rigorouscone} below, Theorem~\ref{thm:generaltwist}
is the statement that $I^\tw(\tau,-z)$ is a
$\Lambda[\![s_0,s_1,\dots]\!][\![\tau]\!]$-valued point of $\cLs$.

\begin{proof}[Proof of Theorem~\ref{thm:generaltwist}] We will assume
  throughout the proof that $F$ is a line bundle and omit the index
  ``${}^{(j)}$'', writing $\rho$ for $\rho^{(j)}$; $f_i$ for
  $f_i^{(j)}$; $N_\theta$ for $N^{(j)}_\theta$; and so on.  The proof
  in the case where $F$ is a direct sum of line bundles requires only
  notational changes.

  Define an element $G_y(x,z)$ in $\CC[y,x,z,z^{-1}][\![s_0,s_1,s_2,\ldots]\!]$ by
  \[ 
  G_{y}(x,z) := \sum_{l,m\ge 0} s_{l+m-1} \frac{B_m(y)}{m!}
  \frac{x^{l}}{l!} z^{m-1}.
  \]
  This satisfies functional equations of gamma-function type:
  \begin{align}
    \label{eq:Gamma1} G_{y}(x,z) &= G_{0}(x+y z, z), \\
    \label{eq:Gamma2} G_{0}(x+z,z) &= G_{0}(x,z) + \bs(x).
  \end{align}
  Equality (\ref{eq:Gamma1}) follows from the fact that the
  coefficient of $s_k$ in $G_y(x,z)$ is the degree $k$ part of
  \[ 
  \left(\sum_{m=0}^\infty \frac{B_m(y)}{m!} z^{m-1} \right)
  \left(\sum_{l=0}^\infty \frac{x^l}{l!}\right ) =
  \frac{e^{x+yz}}{e^z-1},
  \] 
  where $\deg x=\deg z=1$ and $\deg y=0$.  Equality (\ref{eq:Gamma2})
  follows from 
  \[ 
  \frac{e^{x+z}}{e^z-1} = \frac{e^{x}}{e^z-1} + e^x.
  \]

  We need to show that $I^\tw(\tau,-z) \in \cLs$.  As before, write
  elements $\theta \in \NETTX$ as $\theta = (g,d,S)$ with $S =
  (i_1,\ldots,i_n)$.  Observe that
  \begin{align*}
    M_\theta(-z) &= \exp\(
    \sum_{m=K}^{m=N_\theta}
    \bs\(\rho+\(f_{\ibar_n}-m\)z\) 
    -
    \sum_{m=K}^{m=0}
    \bs\(\rho+\(f_{\ibar_n}-m\)z\) 
    \)
    && \text{for $K \ll 0$} \\
    &= \exp \Big(
    G_0\(\rho+f_{\ibar_n}z,z\)
    -
    G_0\(\rho+(f_{\ibar_n}-N_\theta)z,z\)    
    \Big)
    && \text{by \eqref{eq:Gamma2}} \\
   &= \exp \Big(
    G_{f_{\ibar_n}}\(\rho,z\)
    -
    G_0\(\rho+(f_{\ibar_n}-N_\theta)z,z\)    
    \Big)
    && \text{by \eqref{eq:Gamma1}.}
  \end{align*}
  We know that
  \[
  \Delta^\tw = \bigoplus_{i \in \cI} \exp\Big(G_{f_i}(\rho,z)\Big),
  \]
  that $\Delta^\tw\(\cL^\untw\) = \cLs$, and that $J_\theta(\tau,-z)$
  is supported on the component $\cX_{\ibar_n}$ of $\cIX$.  It
  therefore suffices to show that the family
  \[
  \tau \mapsto \sum_{\theta \in \NETTX} \exp \Big(
  - G_0\(\rho+(f_{\ibar_n}-N_\theta)z,z\)    
  \Big) J_\theta(\tau,-z)
  \]
  of elements of $\cH$ lies in the Lagrangian submanifold
  $\cL^\untw$ for the \emph{untwisted} theory.  But
  \begin{multline*}
    \exp \Big(
    - G_0\(\rho+(f_{\ibar_n}-N_\theta)z,z\)    
    \Big) J_\theta(\tau,-z) \\
    =
    \exp \Big(
    - G_0\(\rho - \<\rho,d\>z + f_{i_1} z +f_{i_2} z + \cdots +
    f_{i_{n-1}} z,z\)    
    \Big) J_\theta(\tau,-z),
  \end{multline*}
  and Lemma~\ref{lem:propertiesofJ} shows that this is
  \[
  \exp\Big( - G_0\(z \nabla_\rho + z D,z\) \Big) J_\theta(\tau,-z)
  \]
  where $D = \sum_{i \in \cI} f_i D_i$.  We thus want to show that
  the family
  \begin{equation}
    \label{eq:easyfamily}
    \tau \mapsto \exp\Big( - G_0\(z \nabla_\rho + z D,z\) \Big)
    J^\untw(\tau,-z)
  \end{equation}
  lies on $\cL^\untw$.

  This last statement follows from the geometric properties of
  $\cL^\untw$ established in Appendix~\ref{sec:rigorouscone}.  Let $h$
  be a general point in a formal neighbourhood of $-z$ in $\cH$:
  \begin{align*}
    h=-z+ \sum_{k=0}^\infty t_k z^k+ \sum_{k=0}^\infty \frac{p_k}{(-z)^{k+1}}, &&
    t_k,p_k\in \HorbXL. 
  \end{align*}
  Then $\cL^\untw$ is defined by the equations $E_0=0$, $E_1=0$,
  $E_2=0$,\ldots (\emph{c.f.}  \eqref{eq:pointonL} above and
  (\ref{eq:def_rigorouscone}) below), where
  \[
  E_j(h):= p_j - \sum_{n\ge 0} \sum_{d\in \EffX} \sum_{1\le \epsilon\le
    N} \frac{Q^d}{n!}
  \correlator{\bt(\psi),\dots,\bt(\psi),\psi^j\phi_\epsilon}
  _{0,n+1,d}^{\cX,\untw}
  \phi^\epsilon.
  \]
  Let $\tau \mapsto J_\bs(\tau,-z)$ be the family from
  \eqref{eq:easyfamily}.  The application of $E_j$ to $J_\bs(\tau,-z)$
  is $(\tau,\bs,Q)$-adically convergent; we want to show that it is
  zero.  It is obvious that $E_j(J_\bs(\tau,-z))$ is zero at
  $s_0=s_1=\cdots=0$.  Set $\deg s_k=k+1$ and assume by induction that
  $E_j(J_\bs(\tau,-z))$ vanishes up to degree $n$ in variables
  $s_0,s_1,s_2,\dots$.  Then we have
  \begin{equation}
  \label{eq:diffbys}
  \parfrac{}{s_i} E_j(J_\bs(\tau,-z) ) 
  = d_{J_\bs(\tau,-z)}E_j
  \(z^{-1} P_i(z\nabla,z) J_{\bs}(\tau,-z) \),
  \end{equation} 
  where 
  \[
  P_i(z\nabla,z)=\sum_{m=0}^{i+1} \frac{1}{m! (i+1-m)!}
  z^{m}B_m(0)\(z\nabla_\rho+zD\)^{i+1-m}.
  \] 
  The induction hypothesis shows that we can find a family $\tau
  \mapsto \widetilde{J}_\bs(\tau,-z)$ of elements of $\cL^\untw$
  (\emph{i.e.}, in the language of Appendix~\ref{sec:rigorouscone}, a
  $\Lambda[\![s_0,s_1,\dots]\!][\![\tau]\!]$-valued point of
  $\cL^\untw$) such that $[\widetilde{J}_\bs]_+=[J_\bs]_+$ and that
  $\widetilde{J}_\bs-J_\bs$ consists of terms of degree greater than
  $n$ in $s_0,s_1,s_2,\dots$ ($[J]_+$ discards all negative powers of
  $z$ in $J$).  Then the right hand side of \eqref{eq:diffbys}
  coincides with
  \[
  d_{\widetilde{J}_\bs(\tau,-z)}E_j
  \(z^{-1} P_i(z\nabla,z) \widetilde{J}_{\bs}(\tau,-z) \)
  \]
  up to degree $n$.  But this is zero, as repeated applications of 
  Lemma \ref{lem:diff_arc} and Corollary \ref{cor:Dmod} 
  show that the term in parenthesis is an element of 
  $T_{\widetilde{J}_\bs(\tau,-z)}\cL^\untw$. 
  Hence the left hand side
  of \eqref{eq:diffbys} vanishes up to degree $n$.  This completes the
  induction step, and the proof.
\end{proof}

\section{Application 1: Genus-Zero Invariants of Hypersurfaces}
\label{sec:hypersurfaces}
It is well-known that, for a complete intersection $Y$ which is cut
out of a projective variety $X$ by a section of a direct sum $E \to X$
of convex line bundles, many genus-zero one-point Gromov--Witten
invariants of $Y$ can be obtained as the non-equivariant limit of
twisted genus-zero Gromov--Witten invariants of $X$: one takes $F = E$
and $\bc$ to be the $T$-equivariant Euler class.  This idea lies at
the heart of most proofs of mirror theorems for toric complete
intersections \cite{Givental:equivariant, Kim, Lian--Liu--Yau:1,
  Bertram, Lee, Lian--Liu--Yau:2, Givental:toric, Lian--Liu--Yau:3}.
The same thing holds for complete intersections in orbifolds: given a
direct sum $E \to \cX$ of convex line bundles, the non-equivariant
limit of a genus-zero twisted one-point Gromov--Witten invariant of
$\cX$ (with $F = E$ and $\bc$ the $T$-equivariant Euler class) is a
genus-zero Gromov--Witten invariant of the complete intersection $\cY
\subset \cX$ cut out by a section of $E$.  This is explained in
\cite[Section 5.2]{Tseng}.

\subsection{An Example: A Quintic Hypersurface}
\label{sec:quintic}

We illustrate this in the case of the quintic hypersurface in
$\PP(1,1,1,1,2)$, taking $\cX = \PP(1,1,1,1,2)$; $F$ to be the line
bundle $\cO(5) \to \cX$; $i:\cY \to \cX$ to be the inclusion of the
corresponding hypersurface; and $\bc$ to be the $T$-equivariant Euler
class\footnote{We established notation for $T$-equivariant
  characteristic classes in Section~\ref{sec:twisted}.}.  Let $\be$
denote the non-equivariant Euler class.

The inertia stack $\cIX$ has two components,
\begin{align*}
  &\cX_0 \cong \PP(1,1,1,1,2) && \text{age $0$,} \\
  &\cX_{1 \over 2} \cong \PP(2) && \text{age $2$.}
\end{align*}
If $\fun_i$ is the fundamental class of $\cX_i$  and $p = c_1(\cO(1))$,
then
\[
\phi_0 = \fun_0, \quad
\phi_1 = p \fun_0, \quad
\phi_2 = p^2 \fun_0, \quad
\phi_3 = p^3 \fun_0, \quad
\phi_4 = p^4 \fun_0, \quad
\phi_5 = \fun_{1 \over 2}
\]
is a basis for $\HorbX$.  Theorem~1.6 in \cite{CCLT} shows that the
restriction to the locus $\tau = t p$ of the untwisted $J$-function
$J^\untw(\tau,z)$ of $\cX$ is\footnote{Here and henceforth we write
  $\fr{r}$ for the fractional part of a rational number $r$.}
\[
J^\untw(t p,z) = z e^{t p/z} \sum_{\substack{d:d \geq 0\\2 d \in \ZZ}}
{Q^d e^{d t} \over
\prod_{\substack{b:0<b\leq d\\\fr{b}=\fr{d}}} (p+b z)^4 
\prod_{\substack{b:0<b\leq 2d\\\fr{b}=0}} (2 p+b z) }
\fun_{\fr{d}}.
\]
This is the sum of the contributions $J_\theta(tp,z)$ where the
topological type $\theta = (0,d,S)$ has either $S = (0,0,\ldots,0,0)$,
in which case
\[
M_\theta(z) = \prod_{1 \leq m \leq 5 d} (\lambda_1 + 5 p + m z),
\]
or $S = (0,0,\ldots,0,{1 \over 2})$, in which case 
\[
M_\theta(z) = \prod_{1 \leq m \leq 5 d+{1 \over 2}} (\lambda_1 + 5 p +
(m-\textstyle {1 \over 2} z)).
\]
Thus Theorem~\ref{thm:generaltwist} implies that the family $t \mapsto
I^\tw(tp,-z)$ lies on $\cLs$, where
\[
I^\tw(tp,z) = z e^{t p/z} \sum_{\substack{d:d \geq 0\\2 d \in \ZZ}}
Q^d e^{d t} 
{\prod_{\substack{b: 0 < b \leq 5 d \\ \fr{b} = \fr{d}}} (\lambda_1 + 5 p + b z)
\over
\prod_{\substack{b:0<b\leq d\\\fr{b}=\fr{d}}} (p+b z)^4 
\prod_{\substack{b:0<b\leq 2d\\\fr{b}=0}} (2 p+b z) }
\fun_{\fr{d}}.
\]

We have
\begin{align*}
  I^\tw(tp,z)  &= z e^{t p/z} \( \fun_0 + {60 Q e^t \over z} \fun_0 +
  O(z^{-2}) \) \\
  &= z + tp + 60 Q e^t \fun_0 + O(z^{-1}),
\end{align*}
and it follows, as the twisted $J$-function is characterized by
\eqref{eq:Jcharact}, that
\[
I^\tw(tp,z) = J^\tw(tp + 60 Q e^t \fun_0,z).
\]
But the String Equation \cite[Theorem~8.3.1]{AGV:2} implies that
\[
J^\tw(\tau + a \fun_0,z) = e^{a/z} J^\tw(\tau,z),
\]
and so
\begin{align}
J^\tw(tp,z) &= \exp(-60 Q e^t/z)  I^\tw(tp,z) \nonumber \\
& = 
  z e^{t p/z} \( \fun_0 + {30 Q^{1/2} e^{t/2} \over
    z^2}\fun_{1 \over 2} 
  + {(137 \lambda_1+265 p) Q e^t \over z^2} \fun_0 
  \right. \label{eq:explicitJ}\\
& \qquad \qquad \qquad \qquad \left. + {7650 Q^2 e^{2 t}
  \over z^2} \fun_0 + O(z^{-3}) \) \nonumber .
\end{align}
On the other hand, the Divisor Equation gives
\[
J^\tw(tp,z) = z e^{tp/z}
\( \fun_0 + \sum_{d>0} Q^d e^{d t} 
\<{\phi_\epsilon \over z(z-\psi) }\>^{\cX,\tw}_{0,1,d} \phi^\epsilon \)
\]
where
\[
\phi^i =
\begin{cases}
  {2 p^{4-i} \over \lambda_1 +5 p} \fun_0 & 0 \leq i \leq 4 \\
  2 \fun_{1 \over 2} & i=5
\end{cases}
\]
is the basis for $\HorbX \otimes \CC(\lambda_1)$ which is dual to
$\{\phi_\epsilon\}$ under the twisted pairing
\eqref{eq:twistedpairing}.  Expanding \eqref{eq:explicitJ} in terms of
the $\{\phi^\epsilon\}$, we find that
\begin{align*}
  \< \phi_5 \>^{\cX,\tw}_{0,1,{1 \over 2}} &= 15 & 
  \< \phi_2 \>^{\cX,\tw}_{0,1,1} &= {1325 \over 2} & 
  \< \phi_3 \>^{\cX,\tw}_{0,1,1} &= 475 \lambda_1 \\
  \< \phi_4 \>^{\cX,\tw}_{0,1,1} &= {137 \lambda_1^2 \over 2} &
  \< \phi_3 \>^{\cX,\tw}_{0,1,2} &= 19125 &
  \< \phi_4 \>^{\cX,\tw}_{0,1,2} &= 3825 \lambda_1 .
\end{align*}

We now take the non-equivariant limit $\lambda_1 \to 0$.  Since $F$ is
convex, $F_{0,1,d}$ is a vector bundle.  In the non-equivariant limit,
the twisted Gromov--Witten invariant
\[
\< \phi_\epsilon \>^{\cX,\tw}_{0,1,d} = \( \ev_1^\star \phi_\epsilon \cup
\bc\(F_{0,1,d}\) \) \cap \left[ \cX_{0,1,d} \right]^{\text{\rm vir}}
\]
becomes 
\[
\( \ev_1^\star \phi_\epsilon \cup
\be\(F_{0,1,d}\) \) \cap \left[ \cX_{0,1,d} \right]^{\text{\rm vir}}.
\]
Functoriality for the virtual fundamental class
\cite{Kim--Kresch--Pantev} implies that
\[
\be\(F_{0,1,d}\) \cap \left[ \cX_{0,1,d} \right]^{\text{\rm vir}} = 
j_\star \left[ \cY_{0,1,d} \right]^{\text{\rm vir}}
\]
where $j:\cY_{0,1,d} \to \cX_{0,1,d}$ is the inclusion induced by
$i:\cY \to \cX$, so in the non-equivariant limit
\[
\< \phi_\epsilon \>^{\cX,\tw}_{0,1,d} \longrightarrow \<  i^\star \phi_\epsilon \>^{\cY,\untw}_{0,1,d}.
\]
We conclude, for example, that the virtual number of degree-${1 \over
  2}$ rational curves on $\cY$ --- any such curve passes through the
stacky point on $\cY$ --- is $15$, and that the virtual number of
degree-$1$ rational curves on $\cY$ which meet the cycle dual to $P^2$
is ${1325 \over 2}$.  

\subsection{A Quantum Lefschetz Theorem for Orbifolds}

Let us return now to the general situation of Example~A, so that $F
\to \cX$ is a direct sum of convex line bundles and $\bc$ is the
$T$-equivariant Euler class.  The key point in our analysis of the
quintic hypersurface was that
\[
I^\tw(t,z) = z + f(t) + O(z^{-1}).
\]
This implied the equality $J^\tw(f(t),z) = I^\tw(t,z)$, which
expresses genus-zero twisted invariants (on the left-hand side) in
terms of untwisted invariants (on the right).  In fact, any time we
have an equality of the form
\begin{equation}
  \label{eq:Jscale}
  I^\tw(t,z) = F(t) z \fun_0 + G(t) + O(z^{-1}),
\end{equation}
where $F$ is an invertible scalar-valued function, $G$ takes values in
orbifold cohomology, and $\fun_0 \in \HorbX$ is the identity element,
we can deduce that
\begin{align*}
  J^\tw(\tau(t),z) = {I^\tw(t,z) \over F(t)} && \text{where $\tau(t) =
    {G(t) \over F(t)}$}.
\end{align*}
This follows from the fact that $\cL^\tw$ is (the germ of) a cone
\cite{Tseng,Givental:symplectic} --- so we can divide
\eqref{eq:Jscale} by $F(t)$ and still obtain a family of elements of
$\cL^\tw$ --- and the characterization \eqref{eq:Jcharact} of the
twisted $J$-function.  We can assure an equality \eqref{eq:Jscale},
provided that $c_1(F)$ is not too positive, by restricting $t$ to lie
in an appropriate subspace.

\begin{cor}[Quantum Lefschetz for Orbifolds] \label{thm:orbQL} Suppose
  that $F \to \cX$ is a direct sum of line bundles $F^{(1)}, \ldots,
  F^{(r)}$ such that each of $c_1(F^{(1)})$,\ldots,$c_1(F^{(r)})$, and
  $c_1(\cX) - c_1(F)$ are nef.  Use notation as in
  Section~\ref{sec:QRR} and write $t'$ for a general point in the
  subspace
  \begin{equation}
    \label{eq:goodsubspace}
    \left\{ \alpha \in 
      \bigoplus_{i:f^{(j)}_i =0 \, \forall j} H^\bullet(\cX_i; \CC) : \deg(\alpha) \leq 2
      \right\}
  \end{equation}
  of $\HorbX$.  Then
  \begin{align}
  I^\tw(t',z) & = \sum_{\theta \in \NETTX} \prod_{j=1}^r
  \prod_{m=1}^{N_\theta^{(j)}} \(\lambda_j + \rho^{(j)} + \(m -
    f^{(j)}_{\ibar_n}\)z\)
    \cdot J_\theta(t',z) \label{eq:lefschetzI} \\
    & = F(t') z \fun_0 + G(t') + O(z^{-1}), \label{eq:lefschetzexpansion}
  \end{align}
  for some $F$ and $G$ with $F$ scalar-valued and invertible, and
  \begin{align}
    J^\tw(\tau(t'),z) = {I^\tw(t',z) \over F(t')} && \text{where $\tau(t') =
      {G(t') \over F(t')}$}. \label{eq:lefschetzpunchline}
  \end{align}
\end{cor}

\begin{proof}
  Since $t'$ is supported on those components $\cX_i$ of the inertia
  stack such that each $f_i^{(j)}$ is zero, $J_\theta(t',z)$ vanishes
  unless each $N_\theta^{(j)} \geq 0$.  This proves
  \eqref{eq:lefschetzI}.  The expansion \eqref{eq:lefschetzexpansion}
  follows by computing the highest powers of $z$ that occur in
  $J_\theta(t',z)$ and in the modification factor, and using the
  formula \cite[Theorem~A]{Chen--Ruan:GW} for the virtual dimension of
  the moduli space of stable maps. The rest was explained above.
\end{proof}

If $F$ is in addition convex then we can pass to the non-equivariant
limit, exactly as in Section~\ref{sec:quintic}, and thereby express
genus-zero one-point Gromov--Witten invariants of a complete
intersection $\cY$ cut out by a section of $F$ in terms of the
ordinary Gromov--Witten invariants of $\cX$.  This approach is used in
\cite{CCLT} to compute genus-zero invariants of weighted projective
complete intersections.  If we assume more --- that $H^1(\cC,f^\star
F) = 0$ for all topological types $\theta$ which contribute
non-trivially to $I^\tw(t',z)$ --- then exactly the same argument
allows us to determine those genus-zero $(n+1)$-point Gromov--Witten
invariants of $\cY$ which involve $n$ classes coming from
\eqref{eq:goodsubspace}.

\section{Application 2: Genus-Zero Local Invariants}
\label{sec:local}

Let $G$ be a finite cyclic group, let $\cX = BG$, and let $E \to \cX$
be the vector bundle arising from a representation $\rho:G \to
(\Cstar)^m \subset GL_m(\C)$.  Let $T = (\Cstar)^m$.  The total space
of $E$ is the orbifold $\left[\CC^m / G\right]$, where $G$ acts via
$\rho$, and so the $T$-equivariant Gromov--Witten invariants of
$\left[\CC^m / G\right]$
coincide\footnote{\label{footnote:noncompact}One needs to choose a
  definition of Gromov--Witten invariants of the non-compact orbifold
  $\left[\CC^m / G\right]$.  In the Introduction we defined these to
  be local Gromov--Witten invariants.  One could also define them via
  virtual localization to $T$-fixed points using
  \cite{Graber--Pandharipande}; this gives the same results. } both
with the $T$-equivariant local Gromov--Witten invariants of $E$ and
with twisted Gromov--Witten invariants of $\cX$ where $F=E$ and $\bc$
is the $T$-equivariant inverse Euler class.  In this Section we use
Theorem~\ref{thm:generaltwist} to compute $T$-equivariant
Gromov--Witten invariants of $\left[\CC^2/\ZZ_n\right]$, where $\ZZ_n$
acts with weights $(n-1,1)$, and of $\left[\CC^3/\ZZ_3\right]$ where
$\ZZ_3$ acts with weights $(1,1,1)$.  Our starting point is the
untwisted $J$-function of $B \ZZ_n$.

\subsection{The Untwisted $J$-Function of $B\ZZ_n$} 
\label{sec:JBZ}
Let $\cX = B\ZZ_n$.  Components of the inertia stack $\cIX$ are
indexed by elements of $\ZZ_n$, and hence by the set of fractions
\[
\cI = \left \{ \textstyle{i \over n} : 0 \leq i < n \right\}
\]
via ${r \over n} \in \cI \mapsto [r] \in \ZZ_n$.  Each component of
$\cIX$ is a copy of $B \ZZ_n$, and we write $\fun_i$ for the
fundamental class of the component $\cX_i$.

\begin{pro}
  Let $x = x_0 \fun_0 + x_1 \fun_{1 \over n} + \cdots + x_{n-1}
  \fun_{n-1 \over n} \in H^\bullet_{\text{\rm orb}}(B \ZZ_n;\CC)$. Then
  \[
  J^\untw(x,z)=z \sum_{k_0,k_1,...,k_{n-1}\geq
    0}\frac{1}{z^{k_0 + k_1 + \cdots + k_{n-1}}}\frac{x_0^{k_0} x_1^{k_1}\cdots
  x_{n-1}^{k_{n-1}}}{k_0! \, k_1! \cdots
  k_{n-1}!} \, \fun_{\<\sum_{i=0}^{n-1} {i k_i \over n} \>}.
  \]
\end{pro}

\begin{proof}
  The Gromov--Witten theory of $BG$ for any finite group $G$ has been
  completely solved by Jarvis--Kimura \cite{JK}.  They show in
  particular that the untwisted quantum orbifold product
  $\bullet_\tau$ on $H^\bullet_{\text{\rm orb}}(B \ZZ_n;\CC)$ is
  semisimple and independent of $\tau$.  The untwisted $J$-function is
  the unique solution to the differential equations \eqref{eq:qdes}
  which has the form \eqref{eq:Jcharact}.  Thus
  \[
  J^\untw(\tau,z) = \sum_{\alpha=0}^{n-1}ze^{u^\alpha(\tau)/z}
  \f_\alpha
  \]
  where $\tau = u^0(\tau) \f_0 + \cdots + u^{n-1}(\tau)\f_{n-1}$ is the
  expansion of $\tau$ in terms of the basis of idempotents
  $\{\f_\alpha\}$ for $\bullet_\tau$.  Applying Jarvis--Kimura's
  formula \cite[Proposition~4.1]{JK} for the idempotents completes the
  proof.
\end{proof}

\subsection{Genus-Zero Gromov--Witten Invariants of
  $\left[\CC^2/\ZZ_n\right]$}
\label{sec:An}

Consider now the situation described at the beginning of
Section~\ref{sec:local} in the case where $G = \ZZ_n$, $\cX = B G$,
$m=2$, and $\rho:G \to GL_2(\CC)$ is the representation with weights
$(n-1,1)$.  Twisted Gromov--Witten invariants of $\cX$ here are
$T$-equivariant Gromov--Witten invariants of the type $A$ surface
singularity $\left[\CC^2/\ZZ_n\right]$.

To apply Theorem~\ref{thm:generaltwist}, we need to calculate $I^\tw$.
For $\bk= (k_0,k_1,\dots,k_{n-1}) \in \ZZ^n$, let
\begin{align*}
a(\bk)=  \sum_{i=1}^{n-1}\frac{n-i}{n} k_i  && \text{and} &&
b(\bk) = \sum_{i=1}^{n-1} \frac{i}{n} k_i.
\end{align*}
The term
\[
\frac{1}{z^{k_0 + k_1 + \cdots + k_{n-1}}}\frac{x_0^{k_0} x_1^{k_1}\cdots
  x_{n-1}^{k_{n-1}}}{k_0! \, k_1! \cdots
  k_{n-1}!}\fun_{\<b(\bk)\>}
\]
in the untwisted $J$-function of $\cX$ contributes to $J_\theta(t,z)$
where the topological type $\theta = (0,0,S)$ has $S$ consisting of
some permutation of
\[
\overbrace{0,0,\ldots,0}^{k_0},
\overbrace{\textstyle{1 \over n},{1 \over n},\ldots,{1 \over n}}^{k_1},
\ldots,
\overbrace{\textstyle{n-1 \over n},{n-1 \over n},\ldots,{n-1 \over n}}^{k_{n-1}}
\]
followed by $\<-b(\bk)\>$.  The corresponding modification factor is
\begin{equation*}
M_{k_0,k_1,\ldots,k_{n-1}}(z) :=
\prod_{l=0}^{\fl{a(\bk)}-1}
\(\lambda_1-\(\<a(\bk)\>+l \)z\) 
\prod_{m=0}^{\fl{b(\bk)}-1} 
\(\lambda_2-\(\<b(\bk)\>+m\)z\).
\end{equation*}
Theorem~\ref{thm:generaltwist} implies that the family $x \mapsto
I^\tw(x,-z)$ lies on the Lagrangian submanifold $\cL^\tw$, where $x =
x_0 \fun_0 + x_1 \fun_{1 \over n} + \cdots + x_{n-1} \fun_{n-1 \over
  n}$ and
\begin{equation}
  \label{eq:ItwC2Zn}
  I^\tw(x,z)=z \sum_{k_0,k_1,...,k_{n-1}\geq 0}
  \frac{M_{k_0,k_1,\ldots,k_{n-1}}(z)}{z^{k_0 + k_1 + \cdots + k_{n-1}}}
  \frac{x_0^{k_0}x_1^{k_1}\cdots x_{n-1}^{k_{n-1}}}{k_0! \, k_1!
    \cdots k_{n-1}!}\fun_{\<b(\bk)\>}.
\end{equation}
Since 
\[
\fl{a(\bk)}+\fl{b(\bk)}
=
\begin{cases}
\sum_{i=1}^{n-1}k_i & \text{if $n$ divides $\sum_{i=1}^{n-1}ik_i$} \\
\sum_{i=1}^{n-1}k_i-1  & \text{otherwise}
\end{cases}
\]
we see that
\begin{equation*}
\frac{M_{k_0,k_1,\ldots,k_{n-1}}(z)}{z^{k_1 + \cdots + k_{n-1}}}
= 
{\Gamma\(1 - \<a(\bk)\>\) \over
  \Gamma\(1 - a(\bk) \)} 
{\Gamma\(1-\<b(\bk)\>\) \over
  \Gamma\(1 - b(\bk)\)} z^{-1} 
 + O(z^{-2})
\end{equation*}
unless $n$ divides $\sum_{i=1}^{n-1}ik_i$, in which case
\[
\frac{M_{k_0,k_1,\ldots,k_{n-1}}(z)}{z^{k_1 + \cdots + k_{n-1}}}
= O(z^{-2}).
\]
Thus
\[
I^\tw(x,z) = z + \tau^0 \fun_0 + \tau^1 \fun_{1 \over n} + \cdots +
\tau^{n-1} \fun_{n-1 \over n} + O(z^{-1}),
\]
where 
\begin{equation}
\label{eq:mirrormapA_n}
\tau^r =
\begin{cases}
  x_0 & r=0 \\
 \displaystyle  \sum_{\substack{k_1,...,k_{n-1}\geq 0 : \\
        \fr{b(\bk)} = {r \over n}} } 
  \frac{x_1^{k_1} x_2^{k_2} \cdots x_{n-1}^{k_{n-1}}}{k_1! \, k_2!
    \cdots k_{n-1}!}
  {\Gamma\(1 - \<a(\bk)\>\) \over
    \Gamma\(1 - a(\bk)\)} 
  {\Gamma\(1-\<b(\bk)\>\) \over
    \Gamma\(1 - b(\bk)\)}
  & r \ne 0.
\end{cases}
\end{equation}
Since the twisted $J$-function gives the unique family of elements of
$\cLs$ satisfying \eqref{eq:Jcharact}, it follows that
\begin{equation}
  \label{eq:punchline}
  I^\tw(x,z) = J^\tw(\tau^0 \fun_0 + \tau^1 \fun_{1 \over n} + \cdots +
  \tau^{n-1} \fun_{n-1 \over n},z).
\end{equation}

To calculate genus-zero twisted Gromov--Witten invariants we need to
determine the twisted $J$-function as a function of $\tau^0$, \ldots,
$\tau^{n-1}$, and so we need to invert the \emph{mirror map}
\begin{equation}
  \label{eq:mirrormap}
  (x_0,x_1,\ldots,x_{n-1}) \mapsto (\tau^0,\tau^1,\ldots,\tau^{n-1}).
\end{equation}
In Appendix~\ref{sec:CRC}, we prove:

\begin{pro}[\emph{cf.} Proposition \ref{pro:flatcoordsarebasis}]
\label{pro:inversemirrorA_n}
The inverse to the mirror map \eqref{eq:mirrormap} is given by
\[
x_{i} =
\begin{cases}
  \tau^0 & i=0 \\
  (-1)^{n-i} e_{n-i}(\kappa_0,\kappa_1,\dots, \kappa_{n-1}) &
  i \ne 0
\end{cases}
\]
where $e_j$ is the $j$th elementary symmetric function, $\zeta =
\exp\(\frac{\pi\sqrt{-1}}{n}\)$, and
\[
\kappa_k(\tau^1,\dots,\tau^{n-1}) = \zeta^{2k+1} \prod_{r=1}^{n-1} 
\exp\(\frac{1}{n} \zeta^{(2k+1)r}  \tau^r\).
\]
\end{pro}

This Proposition together with \eqref{eq:punchline} determines closed
formulas for all genus-zero $T$-equivariant non-descendant
Gromov--Witten invariants of $[\CC^2/\ZZ_n]$.  These invariants are
packaged into a generating function called the 
\emph{Gromov--Witten potential} of 
$[\CC^2/\ZZ_n]$ --- see \emph{e.g.}
\cite{Bryan--Graber}*{Section~1.2} for a definition --- which is equal
to
\[
\cF_0^{[\CC^2/\ZZ_n]}(\tau)= \sum_{m \ge 0} {1 \over m!}
\< \tau, \tau,\ldots,\tau \>_{0,m,0}^{\cX,\tw}
\]
where $\cX = B \ZZ_n$ and $\tau = \tau^0 \fun_0 + \cdots + \tau^{n-1}
\fun_{n-1 \over n} $.

\begin{pro}
\label{pro:potA_n}
\[\cF_0^{[\CC^2/\ZZ_n]} = \frac{\(\tau^0\)^3}{6n\lambda_1\lambda_2} + 
\frac{\tau^0}{2n}\sum_{i=1}^{n-1}\tau^i \tau^{n-i} - 
\lambda_1 G(\tau^1,\dots,\tau^{n-1})  - \lambda_2 G(\tau^{n-1},\dots,\tau^{1})
\]
where the derivatives of $G$ are given by 
\begin{multline*}
\parfrac{G}{\tau^{r}}(\tau^1,\dots,\tau^{n-1}) = 
\sum_{\substack{k_1,\dots,k_{n-1}\ge 0 : \\ 
      \fr{b(\bk)} = { n-r \over n} }} 
\frac{(x_1)^{k_1}\cdots (x_{n-1})^{k_{n-1}}}{n \, k_1! \cdots k_{n-1}!} \\
\times \(\sum_{m=0}^{\fl{a(\bk)}-1} \frac{1}{m+\<a(\bk)\>}\)
\frac{\Gamma(1-\<a(\bk)\>)}{\Gamma(1-a(\bk))} 
\frac{\Gamma(1-\<b(\bk)\>)}{\Gamma(1-b(\bk))}
\end{multline*}
and the relationship between $x_i$ and $\tau^r$ is given in Proposition
\ref{pro:inversemirrorA_n}.
\end{pro} 

\begin{proof}
The terms in the potential which involve $\tau^0$ are determined by
\[
\partial_{\tau^0}\partial_{\tau^i}\partial_{\tau^j} \cF_0^{[\CC^2/\ZZ_n]}=
\(\fun_{\frac{i}{n}},\fun_{\frac{j}{n}}\)^{\tw}.
\]
The others can be extracted from the $z^{-1}$ term in
\eqref{eq:punchline}, using the explicit formula \eqref{eq:ItwC2Zn}
for $I^\tw(t,z)$ and the fact that
\begin{equation*}
  J^\tw(\tau,z) = z + \sum_{i=0}^{n-1} \tau^i \fun_{\frac{i}{n}} + \frac{1}{z} 
  \(n\lambda_1 \lambda_2 \parfrac{\cF_0^{[\CC^2/\ZZ_n]}}{\tau^0} \fun_0 +  
  \sum_{i=1}^{n-1} n \parfrac{\cF_0^{[\CC^2/\ZZ_n]}}{\tau^i} \fun_{\frac{n-i}{n}}
  \) + O(z^{-2}).  
\end{equation*}
This equality follows from \eqref{eq:twistedJfun}, as the bases 
  \begin{align*}
    \fun_0,
    \fun_{\frac{1}{n}},
    \fun_{\frac{2}{n}},
    \fun_{\frac{3}{n}},
     \dots,
     \fun_{\frac{n-1}{n}}
    && \text{and} &&
    n \lambda_1\lambda_2\fun_0,
    n\fun_{\frac{n-1}{n}}, 
    n\fun_{\frac{n-2}{n}}, 
    n\fun_{\frac{n-3}{n}}, 
    \dots, 
    n \fun_{\frac{1}{n}}
  \end{align*}
  for $H^{\bullet}_{{\rm orb},T}([\CC^2/\ZZ_n])$ are dual with respect
  to the twisted pairing \eqref{eq:twistedpairing}.
\end{proof}

Proposition~\ref{pro:potA_n} immediately implies an explicit formula
for the differential of $\cF_0^{[\CC^2/\ZZ_n]}$.  When $n=2$, we can
integrate this, recovering a result of Bryan--Graber
\cite{Bryan--Graber}.

\begin{exa*}[{$n=2$, $[\CC^2/\ZZ_2]$}]  
  Propositions \ref{pro:inversemirrorA_n} and \ref{pro:potA_n} give
  \begin{align*}
    x_1 &= -\kappa_0 -\kappa_1 = -\sqrt{-1}( e^{\sqrt{-1}\tau^1/2} -
    e^{-\sqrt{-1}\tau^1/2})
    = 2 \sin \(\frac{\tau^1}{2}\), \\
    \frac{dG}{d\tau^1} & = \sum_{k=1}^{\infty}
    \frac{(x_1)^{2k+1}}{2^{2k+1}} \frac{((2k-1)!!)^2}{(2k+1)!}
    \sum_{m=0}^{k-1} \frac{1}{m+\frac{1}{2}}.
  \end{align*}
  Using $(\frac{d}{d\tau^1})^2 = (\frac{d}{dx_1})^2 - \frac{1}{4} (x_1
  \frac{d}{dx_1})^2$, we find
  \[
  \(\frac{d}{d\tau^1}\)^3 G = \sum_{k=0}^\infty \(
  \frac{x_1}{4}\)^{2k+1} \frac{(2k)!}{(k!)^2} = \frac{x_1}{4}
  \(1-\frac{(x_1)^2}{4}\)^{-1/2} =
  \frac{1}{2}\tan\(\frac{\tau^1}{2}\).
  \]
\end{exa*}

\subsection{Genus-Zero Gromov--Witten Invariants of
  $\left[\CC^3/\ZZ_3\right]$}
\label{sec:C3Z3}

Consider now the situation described at the beginning of
Section~\ref{sec:local} in the case where $G = \ZZ_3$, $\cX = B G$,
$m=3$, and $\rho:G \to GL_3(\CC)$ is the representation with weights
$(1,1,1)$.  Twisted Gromov--Witten invariants of $\cX$ here are
$T$-equivariant Gromov--Witten invariants of
$\left[\CC^3/\ZZ_3\right]$.  We set
\begin{align*}
\alpha(\bk)= \frac{k_1}{3} + \frac{2 k_2}{3}
&& \text{where $\bk = (k_0,k_1,k_2) \in \ZZ^3$.}
\end{align*}

Theorem \ref{thm:generaltwist} implies that the family $x\mapsto
I^\tw(x,-z)$ lies on the Lagrangian submanifold $\cL^\tw$, where $x =
x_0 \fun_0 + x_1 \fun_{1 \over 3} + x_2 \fun_{2 \over 3} $ and 
\[
I^\tw(t,z)=z \sum_{k_0,k_1,k_2\geq 0}
\frac{\prod_{\substack{b:0 \leq b < \alpha(\bk) \\ \fr{b} =
    \fr{\alpha(\bk)}}} (\lambda_1 - bz)(\lambda_2 -bz)(\lambda_3 - bz)}{z^{k_0 + k_1 +  k_2}} 
\frac{x_0^{k_0} x_1^{k_1} x_2^{k_2}}{k_0! \, k_1!
  \, k_2!}\fun_{\<\alpha(\bk)\>}
\]
To obtain an expansion of the form
\[
I^\tw(x,z) = z + f(x) + O(z^{-1})
\]
we restrict to the locus $x_2=0$, obtaining
\[
I^\tw(x_0 \fun_0 + x_1 \fun_{1/3},z) =
z+\tau^0\fun_0+\tau^1\fun_{\frac{1}{3}}+O(z^{-1})
\]
with
\begin{equation}\label{mirrormap_C3_mod_Z3}
\begin{split}
\tau^0 &= x_0 \\
\tau^1&=\sum_{k\geq 0}\frac{(-1)^{3k}(x_1)^{3k+1}}{(3k+1)!}\left(\frac{\Gamma(k+\frac{1}{3})}{\Gamma(\frac{1}{3})}\right)^3.
\end{split}
\end{equation}
The twisted $J$-function is characterized by \eqref{eq:Jcharact}, so
\begin{equation}
  \label{eq:IJC3Z3}
  I^\tw(x_0 \fun_0 + x_1 \fun_{1 \over 3},z) = 
  J^\tw(\tau^0 \fun_0 + \tau^1 \fun_{1 \over 3},z) .
\end{equation}
The $T$-equivariant genus-zero non-descendant potential of
$\left[\CC^3/\ZZ_3\right]$ is equal to
\[
\cF_0^{[\CC^3/\ZZ_3]}(\tau)= \sum_{m \ge 0} {1 \over m!}
\< \tau, \tau,\ldots,\tau \>_{0,m,0}^{\cX,\tw},
\]
where $\cX = B\ZZ_3$ and $\tau = \tau^0 \fun_0 + \tau^1 \fun_{1 \over
  3} + \tau^2 \fun_{2 \over 3}$.

\begin{pro}\label{c3_mod_z3}
We have
\begin{align*}
  {\partial \cF_0^{[\CC^3/\ZZ_3]} \over \partial \tau^1}\(\tau^0 \fun_0
  + \tau^1 \fun_{1 \over 3} \) &= {1 \over 3} \sum_{j \geq 0}
  (-1)^{3j} {(x_1)^{3j+2} \over (3j+2)!} \({\Gamma\(j+{2 \over 3}\)
    \over \Gamma\({2 \over 3}\)}\)^3\\
  {\partial \cF_0^{[\CC^3/\ZZ_3]} \over \partial \tau^2}\(\tau^0 \fun_0
  + \tau^1 \fun_{1 \over 3} \) &= {\tau^0 \tau^1 \over 3} - {1 \over
    3} \sum_{j \geq 0}  {(x_1)^{3j+1} \over (3j+1)!}
  \sum_{r=0}^{r=j-1} {\lambda_1 + \lambda_2 + \lambda_3 \over r + {1
      \over 3}}
\end{align*}
where $\tau^1$ and $x_1$ are related by (\ref{mirrormap_C3_mod_Z3}).
\end{pro}

\begin{proof}
  Since the bases $\fun_{0}, \fun_{1 \over 3}, \fun_{2 \over 3}$ and
  $3 \fun_{0}, 3 \fun_{2 \over 3}, 3 \fun_{1 \over 3}$ for
  $H^{\bullet}_{{\rm orb},T}([\CC^3/\ZZ_3])$ are dual with respect to
  the twisted pairing \eqref{eq:twistedpairing}, we have
  \begin{equation*}
    J^\tw(\tau,z) = z + \tau
    + \frac{3}{z} {\partial \cF_0^{[\CC^3/\ZZ_3]} \over \partial \tau^2}
    \fun_{1 \over 3}
    + \frac{3}{z} {\partial \cF_0^{[\CC^3/\ZZ_3]} \over \partial \tau^1}
    \fun_{2 \over 3}
    + O(z^{-2}).  
  \end{equation*}
  The result follows by equating coefficients of $z^{-1}$ in
  \eqref{eq:IJC3Z3}.
\end{proof}

We do not know\footnote{A combinatorial formula for the inverse has
  recently been given by Bayer and Cadman \cite{Bayer--Cadman}.} how to
invert the mirror map (\ref{mirrormap_C3_mod_Z3}). But one can still
calculate the first few terms of the series expansion for $x_1$ in
terms of $\tau^1$:
\[
x_1=\tau^1+\textstyle\frac{(\tau^1)^4}{648}-\frac{29(\tau^1)^7}{3674160}+\frac{6607(\tau^1)^{11}}{71425670400}-\ldots.
\]
and hence extract genus-zero orbifold Gromov--Witten invariants of
$\left[\CC^3/\ZZ_3\right]$ one-by-one.  For example, if 
\[
N^{\text{\rm orb}}_{0,k} = \< \fun_{1\over 3},  \fun_{1\over 3},
\ldots,  \fun_{1\over 3} \>^{\cX,\tw}_{0,3k,0} 
\]
then Proposition~\ref{c3_mod_z3} gives
\[
\setlength{\extrarowheight}{0.1cm}
\begin{array}{r||c|c|c|c|c|c}
k & 1 & 2 & 3 & 4 & 5 & 6 \\ \hline
N^{\text{\rm orb}}_{0,k} &
\frac{1}{3} &
-\frac{1}{27} &
\frac{1}{9} &
-\frac{1093}{729} &
\frac{119401}{2187} &
-\frac{27428707}{6561}
\end{array}.
\]
This agrees with the predictions of Aganagic--Bouchard--Klemm
\cite[Section 6]{ABK}.

\appendix
\section{The Crepant Resolution Conjecture for Type $A$ Surface
  Singularities}
\label{sec:CRC}

A long-standing conjecture of Ruan states that if $\cX$ is an orbifold
with coarse moduli space $X$ and $Y \to X$ is a crepant resolution
then the small quantum cohomology of $Y$ becomes isomorphic to the
small quantum cohomology of $\cX$ after analytic continuation in the
quantum parameters followed by specialization of some of the
parameters to roots of unity.  A refinement of this conjecture,
proposed recently by Bryan and Graber \cite{Bryan--Graber}, suggests
that if $\cX$ satisfies a Hard Lefschetz condition on orbifold
cohomology then the Frobenius manifold structures defined by the
quantum cohomology of $\cX$ and of $Y$ coincide after analytic
continuation and specialization of parameters (see also
\cite{CCIT:wallcrossings1} for a Hard Lefschetz condition).  This is a
stronger assertion: that the \emph{big} quantum cohomology of $Y$
coincides with that of $\cX$ after analytic continuation plus
specialization, via a linear isomorphism which preserves the
(orbifold) Poincar\'e pairing.  In this Appendix we prove these
conjectures in the case where $\cX$ is the $A_{n-1}$ surface
singularity $\left[\CC^2/\ZZ_n\right]$ and $Y$ is its crepant
resolution.  In fact we prove a more precise statement,
Theorem~\ref{thm:maintheorem} below, which also identifies an
isomorphism and the roots of unity to which the quantum parameters of
$Y$ are specialized.  We learned this statement from Jim Bryan
\citelist{\cite{Bryan:personal} \cite{Bryan--Graber}*{Conjecture~3.1}}
and Fabio Perroni \citelist{\cite{Perroni:personal}
  \cite{Perroni}*{Conjecture~1.9}}.

Our proof of Theorem~\ref{thm:maintheorem} is based on mirror symmetry
for toric orbifolds.  By mirror symmetry we mean the fact, first
observed by Candelas \emph{et~al.} \cite{COGP}, that one can compute
virtual numbers of rational curves in a manifold or orbifold $\cX$ ---
\emph{i.e.}  certain genus-zero Gromov--Witten invariants of $\cX$ ---
by solving Picard--Fuchs equations.  Following Givental we reinterpret
our results from Section~\ref{sec:An} in these terms, observing that
there is a close relationship between a cohomology-valued generating
function for genus-zero Gromov--Witten invariants, called the
\emph{$J$-function of $\cX$}, and a cohomology-valued solution to the
Picard--Fuchs equations called the \emph{$I$-function of $\cX$}.  A
similar relationship holds for $Y$: this is
Proposition~\ref{pro:mirror} below.  After describing the toric
structures of $\cX$ and $Y$ and fixing notation for cohomology and
quantum cohomology, we explain below how to extract the quantum
products for $\cX$ and $Y$ from the Picard--Fuchs equations.  Once we
understand this, Theorem~\ref{thm:maintheorem} follows easily: the
proof is at the end of the Appendix.

A number of cases of Theorem~\ref{thm:maintheorem} were already known.
Ruan's Crepant Resolution Conjecture was established for surface
singularities of type $A_1$ and $A_2$ by Perroni \cite{Perroni}.
Theorem~\ref{thm:maintheorem} was proved in the $A_1$ case by
Bryan-Graber \cite{Bryan--Graber}, in the $A_2$ case by
Bryan--Graber--Pandharipande \cite{Bryan--Graber--Pandharipande}, and
in the $A_3$ case by Bryan--Jiang \cite{Bryan--Jiang}.  
Davesh Maulik has computed the genus-zero Gromov-Witten 
potential of the type A surface singularity $\cX=[\CC^2/\ZZ_n]$ 
for all $n$ (as well as certain higher-genus
Gromov-Witten invariants of $\cX$) \cite{Maulik} and 
the reduced genus-zero Gromov-Witten potential
of the crepant resolution $Y$ \cite{Maulik:Anresolution}; 
Theorem 1 should follow from this. The
quantum cohomology of the crepant resolutions of type ADE surface singularities
has been computed by Bryan-Gholampour \cite{Bryan--Gholampour}.
Skarke \cite{Skarke} and
Hosono \cite{Hosono} have also studied the $A_n$ case, from a point of
view very similar to ours, as part of their investigations of
homological mirror symmetry.

\subsection*{$\cX$ and $Y$ as Toric Orbifolds}
\label{sec:toric}

\begin{figure}
  \centering
  \subfigure[]
  {
    \label{fanforX}
    \begin{picture}(100,90)(-40,-5)
      \put(24,0){$\textstyle
        (1,0)$}
      \put(24,80){$\textstyle (1,n)$}
      \put(0,0){\vector(1,0){20}}
      \put(0,0){\vector(1,4){20}}
    \end{picture}
  }
  \hspace{1cm} 
  \subfigure[]
  {
    \label{fanforY}
    \begin{picture}(100,90)(-40,-5)
      \put(20,40){\makebox(0,0){$\vdots$}}
      \put(24,0){$\textstyle (1,0)$}
      \put(24,20){$\textstyle (1,1)$}
      \put(24,60){$\textstyle (1,n-1)$}
      \put(24,80){$\textstyle (1,n)$}
      \put(75,0){ray $0$}
      \put(75,20){ray $1$}
      \put(75,60){ray $n-1$}
      \put(75,80){ray $n$}
      \put(0,0){\vector(1,0){20}}
      \put(0,0){\vector(1,1){20}}
      \put(0,0){\vector(1,3){20}}
      \put(0,0){\vector(1,4){20}}
    \end{picture}
  }
  \caption{(a) The fan for $\cX$.   (b) The fan for $Y$.}
\end{figure}
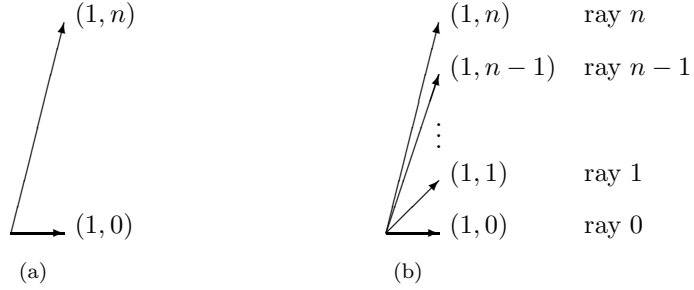

$\cX$ is the toric orbifold corresponding to the fan (or stacky fan
\cite{Borisov--Chen--Smith}) in Figure~\ref{fanforX} and $Y$ is the
toric manifold corresponding to the fan in Figure~\ref{fanforY}.
Background material on toric manifolds and orbifolds can be found in
\cite[Chapter VII]{Audin}.

There is an exact sequence
\[
\begin{CD}
  0 @>>> \ZZ^{n-1} @>{
  M^\mathrm{T}}>>  
  \ZZ^{n+1} @>{
    \begin{pmatrix} \textstyle
      1 & 1 & 1 & \cdots & 1 \\
      0 & 1 &  2 & \cdots & n
    \end{pmatrix}}>> \ZZ^2 @>>> 0,
\end{CD}
\]
and hence we can represent the Gale dual of the right-hand map by
\[
\begin{CD}
  \ZZ^{n+1} @>{M}>> \ZZ^{n-1},
\end{CD}
\]
where 
\[
M= 
 \begin{pmatrix}
   1 & -2 & 1 & 0 & 0 & \cdots & 0 \\
   0 & 1 & -2 & 1  & 0 & \cdots & 0 \\
   \vdots & & \ddots & & \ddots & & \vdots \\
   0 & \cdots & 0 & 1 & -2 & 1 & 0 \\
   0 & \cdots & 0 & 0 & 1 & -2 & 1
 \end{pmatrix}.
\]
Certain faces of the positive orthant $(\RR_{\geq 0})^{n+1} \subset
\RR^{n+1}$ project via $M$ to codimension-$1$ subsets of $\RR^{n-1}$.
The image of the positive orthant is divided by these subsets into
chambers, which are the maximal cones of a fan in $\RR^{n-1}$ called
the \emph{secondary fan} of $Y$.  Chambers in the secondary fan
correspond to toric partial resolutions of $\cX$.  A chamber $K$
corresponds to a fan $\Sigma$ with rays some subset of the rays of the
fan for $Y$, as follows.  Number the rays of the fan for $Y$ as shown
in Figure~\ref{fanforY}.  For a subset $\sigma \subset
\{0,1,\ldots,n\}$, let us write $\bar{\sigma}$ for the complement
$\{0,1,\ldots,n\}\setminus \sigma$, $\RR^\sigma$ for the corresponding
co-ordinate subspace of $\RR^{n+1}$, and say that $\sigma$ covers $K$
iff $K \subset M(\RR^\sigma)$.  The fan $\Sigma$ corresponding to the
chamber $K$ is defined by
\[
  \sigma \in \Sigma  \iff  \text{$\bar{\sigma}$ covers $K$};
\]
the chamber $K$ corresponding to the fan $\Sigma$ is
\[
\bigcap_{\sigma \in \Sigma} M\(\RR^{\bar{\sigma}}\).
\]
We will concentrate on two chambers: $K_{\cX}$, with rays given by the
middle $n-1$ columns of $M$, and $K_Y$ with rays given by the standard
basis vectors for $\RR^{n-1}$.  $K_{\cX}$ corresponds to the toric
orbifold $\cX$ and $K_Y$ corresponds to the toric manifold $Y$.

Let $\cM_{\text{\rm sec}}$ be the toric orbifold corresponding to the
secondary fan of $Y$.  As $K_{\cX}$ and $K_Y$ are simplicial, they
give co-ordinate patches on $\cM_{\text{\rm sec}}$: the co-ordinates
$x_1,\ldots,x_{n-1}$ from $K_{\cX}$ and $y_1,\ldots,y_{n-1}$ from
$K_Y$ are related by
\begin{subequations}
  \label{eq:xtoy}
  \begin{equation}
    y_i =
    \begin{cases}
      x_1^{-2} x_2 & i=1 \\
      x_{i-1} x_i^{-2} x_{i+1} & 1 < i < n-1 \\
      x_{n-2} x_{n-1}^{-2} & i=n-1.
    \end{cases}
  \end{equation}
  More precisely, $x_1,\dots,x_{n-1}$ are multi-valued and the
  co-ordinate patch $\cM_{\rm sec}(K_\cX)$ corresponding to the cone
  $K_{\cX}$ is given by the uniformizing system:
  \[
  \cM_{\rm sec}(K_\cX) \cong \CC^{n-1}/\mu_n, \quad
  (x_1,x_2,\dots,x_{n-1}) \sim (cx_1, c^2 x_2,\dots, c^{n-1}x_{n-1}) \text{
    for } c\in \mu_n.
  \]
  The \emph{$B$-model moduli space $\cM_B$} is the open subset
  $\CC\times \cM_{\rm sec}(K_\cX)$ of $\CC\times \cM_{\rm sec}$.  Denote
  by $x_0$ or $y_0$ the co-ordinate on the first factor $\CC$ of
  $\CC\times \cM_{\rm sec}$, so that
  \begin{equation}
    x_0=y_0.
  \end{equation}
\end{subequations}
We will refer to the point $(x_0,x_1,\ldots,x_{n-1}) = (0,0,\ldots,0)$
as the \emph{large-radius limit point for $\cX$} and the point
$(y_0,y_1,\ldots,y_{n-1}) = (0,0,\ldots,0)$ as the \emph{large-radius
  limit point for $Y$}.  The co-ordinates $x_i$ and $y_j$ are related
to each other by \eqref{eq:xtoy}, so that $y_0,y_1,\ldots,y_{n-1}$ are
co-ordinates on the patch $\CC \times \(\Cstar\)^{n-1} \subset \CC
\times \cM_{\rm sec}(K_\cX) = \cM_B$ where each of
$x_1,x_2,\ldots,x_{n-1}$ is non-zero.

\begin{rem*}
  In what follows the first factor of $\cM_B$, which has co-ordinates
  $x_0$ or $y_0$, will play a rather different role than the second
  factor.  The first factor will correspond under mirror symmetry to
  $H^0_{\text{\rm orb}}(\cX) \subset H^\bullet_{\text{\rm orb}}(\cX)$
  or $H^0(Y) \subset H^\bullet(Y)$, and the second factor will
  correspond to $H^2_{\text{\rm orb}}(\cX) \subset
  H^\bullet_{\text{\rm orb}}(\cX)$ or $H^2(Y) \subset H^\bullet(Y)$.
\end{rem*}

\begin{rem*}
  It would be more honest to define the $B$-model moduli space 
  as the product of $\CC$ with the open subset of $\cM_{\text{\rm
      sec}}$ on which the GKZ system associated to $Y$ is
  non-singular.  This set is slightly smaller than $\cM_B$, as it does
  not contain the discriminant locus of $W_{\cX}$ or $W_Y$ which
  appears below (in the proof of Proposition~\ref{pro:ac}).
\end{rem*}

The presentations of $\cX$ as a toric orbifold and $Y$ as a toric
variety allow us to write $\cX$ and $Y$ as quotients of open sets
$\cU_\cX$, $\cU_Y \subset \CC^{n+1}$ by $\(\Cstar\)^{n-1}$ --- see
\emph{e.g.}  \cite[Chapter VII]{Audin}.  The action of $T =
(\Cstar)^2$ on $\CC^{n+1}$ given by
\begin{equation}
  \label{eq:action}
  (a_0,a_1,\ldots,a_n) \overset{(s,t)}{\longmapsto}
  (sa_0,a_1,a_2,\ldots,a_{n-1}, t a_n)
\end{equation}
descends to give $T$-actions on $\cX$, $X$, and $Y$, and the crepant
resolution $Y \to X$ is $T$-equivariant.  The $T$-fixed locus on $Y$
is the exceptional divisor.  The $T$-action on $\cX =
\left[\CC^2/\ZZ_n\right]$ coincides with that induced by the standard
action of $T$ on $\CC^2$, so the $T$-fixed locus on $\cX$ is the
$B\ZZ_n$ at the origin.  As in the main text, we write
$H_T^\bullet(\{{\rm pt}\}) = \CC[\lambda_1,\lambda_2]$ where
$\lambda_i$ is Poincar\'e-dual to a hyperplane in the $i$th factor of
$\(\mathbf{C P}^\infty\)^ 2 \simeq BT$.

\subsection*{Orbifold Cohomology of $\cX$ and Cohomology of $Y$}
\label{sec:QC}

The $T$-equivariant orbifold cohomology $\HorbTX$ is the
$T$-equivariant cohomology of the inertia stack $\cIX$.  $\cIX$ has
components $\cX_0$, $\cX_1$, \ldots,$\cX_{n-1}$, where
\begin{align*}
  \cX_k = \left[\(\CC^2\)^g / \ZZ_n \right] && \text{with $g = \exp\(2
    k \pi \sqrt{-1} / n\) \in \ZZ_n$.}
\end{align*}
We have
\begin{align*}
  &\cX_k = \left[ \CC^2 / \ZZ_n \right] &  \text{age = $0$} &&&  \text{if $k=0$,} \\
  & \cX_k = B \ZZ_n &  \text{age = $1$} &&&  \text{otherwise.}
\end{align*}
Let $\delta_i$ be the fundamental class of $\cX_i$, $0 \leq i < n$;
this gives a $\CC[\lambda_1,\lambda_2]$-basis for $\HorbTX$.  The
pullback of $\delta_i$ along the inclusion of the $T$-fixed locus
$B\ZZ_n \to \cX$ is the class $\fun_{i \over n}$ from
Section~\ref{sec:JBZ}.  The canonical involution $I$ on $\cIX$ fixes
$\cX_0$ and exchanges $\cX_i$ with $\cX_{n-i}$, $1 \leq i < n$.  As
$I$ is age-preserving, $\HorbX$ satisfies Hard Lefschetz
\citelist{\cite{Bryan--Graber}*{Definition~1.1} \cite{Fernandez}}.

The cone $K_Y$ is the K\"ahler cone for $Y$ and its rays determine a
basis $\gamma_1,\ldots,\gamma_{n-1}$ for $H^2(Y;\ZZ)$.  The dual basis
$\beta_1,\ldots,\beta_{n-1}$ for $H_2(Y;\ZZ)$ is \emph{positive} in
the sense of \cite[Section~1.2]{Bryan--Graber}.  If we define
$\gamma_0 = 1$ and choose lifts of $\gamma_1,\ldots, \gamma_{n-1}$ to
$T$-equivariant cohomology then $\gamma_i$, $0 \leq i < n$, is an
$\CC[\lambda_1,\lambda_2]$-basis for $\HTY$.  We choose a standard
equivariant lift of each $\gamma\in H^2(Y;\ZZ)$ in the following way.
There is a unique representation $\rho_\gamma$ of $\(\Cstar\)^{n-1}$
such that $\gamma$ is the first Chern class of the line bundle
\[
L_\gamma:=\cU_Y\times_{\rho_\gamma}\CC \longrightarrow \cU_Y/(\Cstar)^{n-1} =Y.
\]
This line bundle $L_\gamma$ admits a $T$-action such that $T$ acts on
$\cU_Y$ via \eqref{eq:action} and acts trivially on the $\CC$ factor, and
the lift $\gamma\in H^2_T(Y;\ZZ)$ is the $T$-equivariant first Chern
class of $L_\gamma$.  The columns of $M$, together with the action
\eqref{eq:action}, define elements $\omega_j \in H^2_T(Y;\CC)$, $0
\leq j \leq n$, where
\[
\omega_j =
\begin{cases}
  \lambda_1 + \gamma_1 & j=0 \\
  -2 \gamma_1 + \gamma_2 & j=1 \\
  \gamma_{j-1} - 2 \gamma_j + \gamma_{j+1} & 1<j<n-1 \\
  \gamma_{n-2} - 2 \gamma_{n-1} & j=n-1 \\
  \lambda_2 + \gamma_{n-1} & j=n.
\end{cases}
\]
The class $\omega_i$ is the $T$-equivariant Poincar\'{e} dual of the toric
divisor given in co-ordinates \eqref{eq:action} by $a_i=0$.  We have
\begin{align*}
  \HTY &= \CC[\lambda_1,\lambda_2,\gamma_1,\ldots,\gamma_{n-1}]/ \left
    \langle \omega_i \omega_j : i-j>1 \right \rangle.
\end{align*}

$\cX$ and $Y$ are non-compact but nonetheless one can define
(orbifold) Poincar\'e pairings on the localized $T$-equivariant
(orbifold) cohomology groups
\begin{align*}
  \HorblocTX := H^\bullet_{T,\text{\rm orb}}(\cX;\CC) \otimes
  \CC(\lambda_1,\lambda_2) 
  && \text{and} && 
  \HlocTY := H_T^\bullet(Y;\CC)\otimes \CC(\lambda_1,\lambda_2) 
\end{align*}
using the Bott residue formula.  These pairings take values in
$\CC(\lambda_1,\lambda_2)$, and are non-degenerate.  We write
$\{\gamma^i\}$ and $\{\delta^i\}$ for the bases dual respectively to 
$\{\gamma_i\}$ and $\{\delta_i\}$ under these pairings.

\subsection*{Gromov--Witten Invariants of $\cX$ and $Y$}

As discussed in \cite{Bryan--Graber}, even though some moduli spaces
of stable maps to $\cX$ or $Y$ are non-compact the $T$-fixed loci on
these moduli spaces are compact and so we can still define
$\CC(\lambda_1,\lambda_2)$-valued Gromov--Witten invariants of $\cX$
and $Y$ using the virtual localization formula of
Graber--Pandharipande \cite{Graber--Pandharipande}.  For
$\alpha_1,\ldots,\alpha_n \in \HlocTY$, $d \in H_2(Y;\ZZ)$, and
$i_1,\ldots,i_n \geq 0$, we set
\[
\correlator{\alpha_1 \psi^{i_1}, \ldots, \alpha_n
  \psi^{i_n}}^Y_d = 
\int_{\left[Y_{0,n,d}\right]^{\text{\rm vir}}} \prod_{j=1}^n \ev_j^\star
  \alpha_j \cdot \psi_j^{i_j}.
\]
Here $\psi_i$ and $Y_{0,n,d}$ are as in Section~\ref{sec:twisted} and
the integral is defined by localization to the $T$-fixed substack, as
in \citelist{\cite{Graber--Pandharipande}*{Section~4}
  \cite{Coates}*{Section~3.1}}.  We make a similar definition for
$\cX$; the discussion around footnote \ref{footnote:noncompact} in the
main text explains how to express the resulting correlators
$\correlator{\alpha'_1 \psi^{i_1}, \ldots, \alpha'_n
  \psi^{i_n}}^\cX_0$ as twisted Gromov--Witten invariants of $B\ZZ_n$.

\subsection*{Small Quantum Cohomology and Ruan's Conjecture}

The small quantum product for $\cX$ is the
$\CC(\lambda_1,\lambda_2)$-algebra defined by
\begin{equation}
  \label{eq:smallQCX}
  \delta_i \underset{\rm small}{\star} \delta_j = \sum_{k=0}^{n-1} 
  \correlator{\delta_i,\delta_j,\delta_k}^\cX_{0} 
  \delta^k.
\end{equation}
This coincides with the Chen--Ruan orbifold cup product on
$\HorblocTX$ \cite{Chen--Ruan:orbifold}.  The small quantum product
for $Y$ is the family of $\CC(\lambda_1,\lambda_2)$-algebras,
depending on parameters $q_1,\ldots,q_{n-1}$, defined by
\begin{equation}
  \label{eq:smallQCY}
  \gamma_i \underset{\rm small}{\star} \gamma_j = \sum_{d} \sum_{k=0}^{n-1} 
  \correlator{\gamma_i,\gamma_j,\gamma_k}^Y_{d} 
  q_1^{d_1} \cdots q_{n-1}^{d_{n-1}} \gamma^k.
\end{equation}
where the sum is over classes $d = d_1 \beta_1 + \cdots + d_{n-1}
\beta_{n-1}$ with each $d_i \geq 0$.  It can be obtained from the
untwisted product defined in \eqref{eq:twistedmultiplication} by
setting $\tau=0$ and $Q^d = q_1^{d_1} \cdots q_{n-1}^{d_{n-1}}$.  This
change in notation reflects a change in perspective: in
\eqref{eq:twistedmultiplication} the variable $Q$ was part of the
ground ring $\Lambda$ and we thought of the product $\bullet_\tau$ as
depending formally on $Q$ and $\tau$; here the ground ring is
$\CC(\lambda_1,\lambda_2)$ not $\Lambda$ and we think of
$\underset{\rm small}{\star}$ as a family of products on $\HlocTY$
which varies analytically with $q_1,\ldots,q_{n-1}$.  It follows from
the discussion below that the right-hand side of \eqref{eq:smallQCY}
converges to an analytic function of $q_1,\ldots,q_{n-1}$ in some
neighbourhood of the origin.  Ruan's conjecture asserts that there is
a linear isomorphism $\HorblocTX \to \HlocTY$ which identifies the
products \eqref{eq:smallQCX} and \eqref{eq:smallQCY} after analytic
continuation in the $q_i$ followed by setting the $q_i$ equal to
certain roots of unity.

\subsection*{Big Quantum Cohomology and the Bryan--Graber Conjecture}

The big quantum cohomology of $\cX$ is the family of
$\CC(\lambda_1,\lambda_2)$-algebras parametrized by $\tau \in
\HorblocTX$, $\tau = \tau^0 \delta_0 + \tau^1 \delta_1 + \cdots +
\tau^{n-1} \delta_{n-1}$, defined by
\begin{equation}
  \label{eq:QCX}
  \delta_i \underset{\rm big}{\star} \delta_j = \sum_{m=0}^\infty \sum_{k=0}^{n-1} 
  {1 \over m!} \big\langle
  \delta_i,\delta_j,\delta_k,\overbrace{\tau,\ldots,\tau}^m
  \big\rangle^\cX_{0} 
  \delta^k.
\end{equation}
The big quantum cohomology of $Y$ is the family of
$\CC(\lambda_1,\lambda_2)$-algebras parametrized by $t \in \HlocTY$,
$t = t^0 \gamma_0 + t^1 \gamma_1 + \cdots + t^{n-1} \gamma_{n-1}$,
defined by
\begin{equation}
  \label{eq:QCY}
  \gamma_i \underset{\rm big}{\star}  \gamma_j = 
  \sum_{d} \sum_{m=0}^\infty \sum_{k=0}^{n-1} 
  {1 \over m!} \big\langle
  \gamma_i,\gamma_j,\gamma_k,\overbrace{t,\ldots,t}^m
  \big\rangle^Y_{d} 
  \gamma^k.
\end{equation}
The first sum here is over classes $d = d_1 \beta_1 + \cdots + d_{n-1}
\beta_{n-1}$ with each $d_i \geq 0$.  It follows from the discussion
below that the right-hand sides of \eqref{eq:QCX} and respectively
\eqref{eq:QCY} converge to analytic functions of
$\tau^0,\ldots,\tau^{n-1}$ and respectively $t^0,\ldots,t^{n-1}$ on
appropriate domains.  Note that the product \eqref{eq:QCY} differs
from the untwisted product $\bullet_t$ defined in
\eqref{eq:twistedmultiplication} as it does not contain factors of
$Q^d$.  This is better for our purposes, as the Divisor Equation
implies that \eqref{eq:twistedmultiplication} contains redundant
information: see \cite{Bryan--Graber}*{Section~2.2} for a discussion
of this.

Together with the (orbifold) Poincar\'e pairings, the big quantum
products \eqref{eq:QCX} and \eqref{eq:QCY} define Frobenius
manifolds\footnote{These Frobenius manifolds are defined over the
  field $\CC(\lambda_1,\lambda_2)$.} based on $\HorblocTX$ and
$\HlocTY$.  The Bryan--Graber version of the Crepant Resolution
Conjecture asserts that these Frobenius manifolds coincide after
analytic continuation in the $t^i$ and an appropriate
change-of-variables.  This is our main result.

\begin{thm} \label{thm:maintheorem} The big quantum products
  \eqref{eq:QCX} for $\cX$ and \eqref{eq:QCY} for $Y$ coincide after
  analytic continuation in the $t^i$, the affine-linear change-of-variables
  \[
  t^i =
  \begin{cases}
    \tau^0, & i=0 \\
    -{2 \pi \sqrt{-1} \over n} + \sum_{j=1}^{n-1} L^i_{\phantom{i} j} \tau^j,  & i>0,
  \end{cases}
  \]
  where
  \begin{align*}
    L^i_{\phantom{i} j} = {\zeta^{2 i j} \( \zeta^{-j} - \zeta^{j}\) \over n}, &&  
    \zeta = \exp\({\pi \sqrt{-1} \over n}\),
  \end{align*}
  and the linear isomorphism 
  \begin{equation}
    \label{eq:lineariso}
    \begin{aligned}
      L:\HorblocTX & \to \HlocTY \\
      \delta_0 & \mapsto \gamma_0,\\
      \delta_j & \mapsto \sum_{i=1}^{n-1} L^i_{\phantom{i} j} \gamma_i, && \qquad 1 \leq j < n.
    \end{aligned}
  \end{equation}
  Furthermore, the isomorphism \eqref{eq:lineariso} matches the
  Poincar\'e pairing on $\HlocTY$ with the orbifold Poincar\'e pairing
  on $\HorblocTX$.
\end{thm}

Theorem~\ref{thm:maintheorem} establishes Conjecture~3.1 in
\cite{Bryan--Graber} for the case of polyhedral and binary polyhedral
groups of type $A$, and also Conjecture~1.9 in \cite{Perroni}.  The
path along which analytic continuation is taken is described after
Proposition~\ref{pro:ac} below.

The Divisor Equation implies that we can write \eqref{eq:QCY} as
\[
 \gamma_i \underset{\rm big}{\star}  \gamma_j = \sum_{d} \sum_{k=0}^{n-1} 
  \correlator{\gamma_i,\gamma_j,\gamma_k}^Y_{d} 
  e^{d_1 t_1 + \cdots + d_{n-1} t_{n-1}} \gamma^k.
\]
To pass from the big quantum cohomology algebras of $\cX$ and $Y$ to
the small quantum cohomology algebras, therefore, set $\tau^i = 0$ and
$e^{t^i} = q_i$, $1 \leq i < n$.

\begin{cor}
  The small quantum products \eqref{eq:smallQCX} for $\cX$ and
  \eqref{eq:smallQCY} for $Y$ coincide after analytic continuation in
  the $q_i$, the linear isomorphism \eqref{eq:lineariso}, and the
  specialization 
  \begin{align*}
    q_i = \exp\(-{2 \pi \sqrt{-1} \over n}\), && 1 \leq i < n.
  \end{align*}
\end{cor}

Note that from this point of view the specialization $q_i = c_i$ of
quantum parameters to roots of unity arising in Ruan's conjecture just
reflects the affine-linear identification of flat co-ordinates
\begin{align*}
  t^i = \log c_i + \sum_{j=1}^{n-1} L^i_{\phantom{i} j} \tau^j, && 1 \leq i < n.
\end{align*}

\subsection*{Mirror Symmetry}
\label{sec:mirrorsymmetry}

As discussed above, by mirror symmetry we mean the fact that one can
compute certain genus-zero Gromov--Witten invariants of $\cX$ and $Y$
by solving Picard--Fuchs equations.  We now make this precise.
Following Givental we compare two cohomology-valued generating
functions for genus-zero Gromov--Witten invariants, the
\emph{$J$-functions} of $\cX$ and $Y$, with two cohomology-valued
solutions to the Picard--Fuchs equations called the
\emph{$I$-functions} of $\cX$ and $Y$.  Mirror symmetry for us, in
this situation, is the statement that the $I$-function of $\cX$ (or
$Y$) coincides with the $J$-function of $\cX$ (or $Y$) after a change
of variables.  After proving this, which is
Proposition~\ref{pro:mirror} below, we then describe how to extract
the quantum products \eqref{eq:QCX} and \eqref{eq:QCY} from the
Picard--Fuchs equations, and finally explain how this implies
Theorem~\ref{thm:maintheorem}.

\subsection*{The $J$-Functions of $\cX$ and $Y$}

The $J$-function of $\cX$ is
\[
J_\cX(\tau,z)= 
z + \tau + 
 \sum_{m=0}^\infty \sum_{k=0}^{n-1} 
  {1 \over m!} \Big\langle
  \overbrace{\tau,\ldots,\tau}^m,
  {\delta_k \over z - \psi}
  \Big\rangle^\cX_{0} 
  \delta^k
\]
where we expand ${1 / (z-\psi)}$ as $\sum_m {\psi^m / z^{m+1}}$.
$J_\cX(\tau,z)$ is a function of $\tau \in \HorblocTX$, $\tau = \tau^0
\delta_0 + \cdots + \tau^{n-1} \delta_{n-1}$, which takes values in
$\HorblocTX \otimes \CC(\!(z^{-1})\!)$.  It is defined and analytic in
an open subset of $\HorblocTX$ where $|\tau^1|,\ldots,|\tau^{n-1}|$
are sufficiently small; this follows from Proposition~\ref{pro:mirror}
below.

The $J$-function  of $Y$ is
\[
J_Y(t,z) = 
z + t + 
\sum_{d} \sum_{m=0}^\infty \sum_{k=0}^{n-1} 
{1 \over m!} \Big\langle
\overbrace{t,\ldots,t}^m,
{\gamma_k \over z - \psi}
\Big\rangle^Y_{d} 
\gamma^k
\]
where the first sum is over $d = d_1 \beta_1 + \cdots + d_{n-1}
\beta_{n-1}$ with each $d_i \geq 0$.  $J_Y(t,z)$ is a function of $t
\in \HlocTY$, $t = t^0 \gamma_0 + \cdots +t^{n-1} \gamma_{n-1}$, which
takes values in $\HlocTY \otimes \CC(\!(z^{-1})\!)$. It is defined and
analytic in an open subset of $\HlocTY$ where $\Re (t^i) \ll 0$, $1
\leq i < n$; this again follows from Proposition~\ref{pro:mirror}.

\begin{rem*}
  $J_\cX(\tau,z)$ differs from the twisted $J$-function of $B \ZZ_n$
  defined in Section~\ref{sec:twistedJ} only in that there we regarded
  the twisted $J$-function as a formal series in
  $\tau^0,\ldots,\tau^{n-1}$ and here we regard $J_\cX(\tau,z)$ as an
  analytic function of the $\tau^i$.  $J_Y(t,z)$ differs from the
  untwisted $J$-function of $Y$ defined in Section~\ref{sec:twistedJ}
  in the same way, and also in that $J_Y(t,z)$ contains no factors of
  $Q^d$.  Using the String Equation and the Divisor Equation, we can
  write $J_Y(t,z)$ as
  \[
  e^{t^0/z} e^{\(t^1\gamma_1+ \cdots + t^{n-1} \gamma_{n-1}\)/z} \( z
  \gamma_0 + \sum_{d} \sum_{k=0}^{n-1} \correlator{\gamma_k \over z -
    \psi}^Y_{d} e^{d_1 t^1 +\cdots + d_{n-1} t^{n-1}} \gamma^k \)
  \]
  and so our definition of $J_Y$ agrees, up to a factor of $z$, with
  that in \cite{Givental:toric}*{Section~1}.
\end{rem*}

\subsection*{The $I$-Functions of $\cX$ and $Y$}

For $\bk= (k_1,\dots,k_{n-1}) \in \ZZ^{n-1}$, let
\[
D_j(\bk) =
\begin{cases}
   -{1 \over n} \sum_{i=1}^{n-1} (n-i) k_i & j=0 \\
   k_j & 1 \leq j < n \\
   -{1 \over n} \sum_{i=1}^{n-1} i k_i & j=n.
\end{cases}
\]
Let $i(\bk) = \fr{-D_n(\bk)}$, where $\fr{r}$ denotes the fractional
part of a rational number $r$.  The $I$-function $I_\cX(x,z)$ of $\cX$
is defined to be
\[
z e^{x_0/z}\, \sum_{k_1,\ldots,k_{n-1} \geq 0}
\frac{
\prod_{\substack{r:D_0(\bk) < r \leq 0 \\ \fr{r} = \fr{D_0(\bk)}}} 
(\lambda_1 + r z)
\prod_{\substack{s:D_n(\bk) < s \leq 0 \\ \fr{s} = \fr{D_n(\bk)}}} 
(\lambda_2 + s z)
}
{z^{k_1 + \cdots + k_{n-1}}}
{ x_1^{k_1} \cdots x_{n-1}^{k_{n-1}}
  \over
  k_1! \cdots k_{n-1}! }
\, \delta_{i(\bk)} .
\]
This is a function of $x = (x_0,\ldots,x_{n-1}) \in \cM_B$, $z \in
\Cstar$, and $\lambda_1,\lambda_2 \in \CC$ which takes values in
$\HorbTX$.  Each component of $I_\cX(x,z)$ with respect to the basis
$\{\delta_i\}$ is an analytic function of $(x,z,\lambda_1,\lambda_2)$
defined in a domain where $|x_1|,\dots,|x_n|$ are sufficiently small
and $x_0,z, \lambda_1,\lambda_2$ are arbitrary.  By taking a Laurent
expansion at $z = \infty$ we can regard $I_\cX(x,z)$ as an analytic
function of $(x,\lambda_1,\lambda_2)$ which takes values in
$\HorblocTX\otimes\CC(\!(z^{-1})\!)$.  $I_\cX(x,z)$ satisfies a system
of Picard--Fuchs equations, as follows.  Define differential operators
\[
\beth_j =
\begin{cases}
   \lambda_1-{1 \over n} \sum_{i=1}^{n-1} (n-i) z x_i \parfrac{}{x_i} & j=0 \\
   z x_j\parfrac{}{x_j} & 1 \leq j < n \\
   \lambda_2 -{1 \over n} \sum_{i=1}^{n-1} i z x_i\parfrac{}{x_i} & j=n.
\end{cases}
\]
Then
\begin{subequations}
  \label{eq:PFX}
  \begin{multline}
    \label{eq:PFXa}
    \( \prod_{j:D_j(\bk)>0} \prod_{m=0}^{D_j(\bk)-1} \(\beth_j - m z\) \) I_{\cX}(x,z) = \\
    x_1^{k_1} \cdots x_{n-1}^{k_{n-1}} \(
    \prod_{j:D_j(\bk)<0} \prod_{m=0}^{-D_j(\bk)-1} \(\beth_j - m z\)
    \) I_{\cX}(x,z).
  \end{multline}
  for each $\bk= (k_1,\dots,k_{n-1}) \in \ZZ^{n-1}$ such that $i(\bk)
  = 0$, and
  \begin{equation}
    \label{eq:PFXb}
    z \parfrac{}{x_0} I_{\cX}(x,z) = I_{\cX}(x,z).
  \end{equation}
\end{subequations}

The $I$-function of $Y$ is
\[
I_Y(y,z) = 
z \,e^{y_0/z}  y_1^{\gamma_1/z} \cdots y_{n-1}^{\gamma_{n-1}/z}  
\sum_d 
\prod_{j=0}^{n} \textstyle {\prod_{m \leq 0} (\omega_j + m z) \over 
\prod_{m \leq D'_j(d)} (\omega_j + m z)}
  y_1^{d_1} \cdots y_{n-1}^{d_{n-1}},
\]
where $y_i^{\gamma_i/z} = \exp\(\gamma_i\log y_i/z\)$, the sum is
over $d = d_1 \beta_1 + \cdots + d_{n-1} \beta_{n-1}$ with each
$d_i \geq 0$, and
\[
D'_j(d) = 
\begin{cases}
  d_1 & j=0 \\
  -2 d_1 + d_2 & j=1 \\
  d_{j-1} - 2 d_j + d_{j+1} & 1<j<n-1 \\
  d_{n-2} - 2 d_{n-1} & j=n-1 \\
  d_{n-1} & j=n.
\end{cases}
\] 
$I_Y(y,z)$ is a multi-valued function of $y = (y_0,\ldots,y_{n-1}) \in
\cM_B$, $z \in \Cstar$, and $\lambda_1,\lambda_2 \in \CC$ which takes
values in $\HTY$.  Each component of $I_Y(y,z)$ with respect to the
basis $\{\gamma_i\}$ is a multi-valued analytic function of
$(y,z,\lambda_1,\lambda_2)$ defined in a domain where $|y_1|,\dots,
|y_{n-1}|$ are sufficiently small, $|z| >
\max(|\lambda_1|,|\lambda_2|)$, and $y_0$ is arbitrary.  By taking a
Laurent expansion at $z = \infty$ we can regard $I_Y(y,z)$ as a
multi-valued analytic function of $(y, \lambda_1,\lambda_2)$ which
takes values in $\HlocTY\otimes \CC(\!(z^{-1})\!)$.  It also satisfies
a system of Picard--Fuchs equations.  Define differential operators
\[
\daleth_j =
\begin{cases}
  \lambda_1 + z y_1 \parfrac{}{y_1} & j=0 \\
  -2 z y_1 \parfrac{}{y_1} + z y_2 \parfrac{}{y_2} & j=1 \\
  z y_{j-1} \parfrac{}{y_{j-1}} - 2 z y_j \parfrac{}{y_j} + z y_{j+1} \parfrac{}{y_{j+1}} & 1<j<n-1 \\
  z y_{n-2} \parfrac{}{y_{n-2}} - 2 z y_{n-1} \parfrac{}{y_{n-1}} & j=n-1 \\
  \lambda_2 + z y_{n-1} \parfrac{}{y_{n-1}} & j=n.
\end{cases}
\]
Then 
\begin{subequations}
  \label{eq:PFY}
  \begin{multline}
    \label{eq:PFYa}
    \( \prod_{j:D'_j(d) >0} \prod_{m=0}^{D'_j(d) - 1}
    \(\daleth_j - m z\)
    \) I_Y(y,z) \\
    = y_1^{d_1} \cdots y_{n-1}^{d_{n-1}} \( \prod_{j:D'_j(d) <0}
    \prod_{m=0}^{-D'_j(d) - 1} \(\daleth_j - m z\) \) I_Y(y,z)
  \end{multline}
  for every $d = d_1 \beta_1 + \cdots + d_{n-1} \beta_{n-1}$, and
  \begin{equation}
    \label{eq:PFYb}
    z \parfrac{}{y_0} I_{Y}(y,z) = I_{Y}(y,z).
  \end{equation}
\end{subequations}
The Picard--Fuchs systems \eqref{eq:PFX} for $\cX$ and \eqref{eq:PFY}
for $Y$ coincide under the co-ordinate change \eqref{eq:xtoy}.  Thus
there is a global system of Picard--Fuchs equations --- a $\cD$-module
over all of $\cM_B$ --- which gives \eqref{eq:PFX} near the
large-radius limit point for $\cX$ and \eqref{eq:PFY} near the
large-radius limit point for $Y$.  This global nature of the
Picard--Fuchs system will play a key role in what follows.

\subsection*{A Mirror Theorem}

By mirror symmetry, we mean the following.

\begin{pro}  \label{pro:mirror} \

  \begin{enumerate}
  \item $I_\cX(x,z)$ and $J_\cX(\tau,z)$ coincide after a change of
    variables expressing $\tau$ in terms of $x$.
  \item $I_Y(y,z)$ and $J_Y(t,z)$ coincide after a change of
    variables expressing $t$ in terms of $y$.
  \end{enumerate}
\end{pro}

\begin{proof}
  Part (1) is equation \eqref{eq:punchline}: $I^\tw(x,z)$ there
  coincides with $I_\cX(x,z)$ here and $J^\tw(\tau,z)$ there coincides
  with $J_\cX(\tau,z)$ here.  The argument that proves Theorem~0.2 in
  \cite{Givental:toric} also proves part (2) here.  Theorem~0.2 as
  stated only applies to compact semi-positive toric manifolds, but
  the proof applies essentially without change to the non-compact
  toric Calabi--Yau manifold $Y$.
\end{proof}

\begin{rem*}
  We learned from Bong Lian that, in unpublished work, he and
  Chien-Hao Liu have established mirror theorems for non-compact toric
  Calabi--Yau manifolds using the arguments of
  \cite{Lian--Liu--Yau:3}.  Once again, the proof for compact toric
  manifolds applies also to the non-compact toric Calabi--Yau case
  without significant change.  This gives an alternative proof of the
  second part of Proposition~\ref{pro:mirror}.
\end{rem*}

From equation \eqref{eq:mirrormapA_n} we know that the change of
variables in Proposition~\ref{pro:mirror} between $x_0,\ldots,x_{n-1}$
and $\tau = \tau^0 \delta_0 + \cdots + \tau^{n-1} \delta_{n-1}$ is
$\tau^r = f^r(x)$,
\[
f^r(x) = 
\begin{cases}
  x_0 & r=0 \\
 \displaystyle  \sum_{\substack{k_1,...,k_{n-1}\geq 0 : \\
        \fr{b(\bk)} = {r \over n}} } 
  \frac{x_1^{k_1} x_2^{k_2} \cdots x_{n-1}^{k_{n-1}}}{k_1! \, k_2!
    \cdots k_{n-1}!}
  {\Gamma\(\<D_0(\bk)\>\) \over
    \Gamma\(1 + D_0(\bk)\)} 
  {\Gamma\(\<D_n(\bk)\>\) \over
    \Gamma\(1 + D_n(\bk)\)}
  & r \ne 0.
\end{cases}
\]
As
\[
f^r(x) = x_r + \text{quadratic and higher order terms in
  $x_1,\ldots,x_{n-1}$}
\]
the functions $f^0(x),\ldots,f^{n-1}(x)$ define co-ordinates on a
neighbourhood of the large-radius limit point for $\cX$ in $\cM_B$.
We call these \emph{flat co-ordinates for $\cX$}.  Similarly, 
\[
J_Y(t,z) = z + t^0 \gamma_0 + t^1 \gamma_1 + \cdots + t^{n-1}
\gamma_{n-1} + O(z^{-1})
\]
and
\[
I_Y(y,z) = z +  g^0(y) \gamma_0 + g^1(y) \gamma_1 + \cdots + g^{n-1}(y) \gamma_{n-1} +
O(z^{-1}) 
\]
for some functions $g^0(y), \ldots, g^{n-1}(y)$ with $g^0(y) = y_0$ and for $1 \leq k < n$,
\[
g^k(y) = \log y_k + \text{single-valued analytic function of
  $y_1,\ldots,y_{n-1}$}.
\]
The change of variables which equates $I_Y$ and $J_Y$ is $t^i =
g^i(y)$, $0 \leq i < n$.  The functions $g^0(y), \ldots, g^{n-1}(y)$
define multi-valued co-ordinates on a neighbourhood of the
large-radius limit point for $Y$; these are the \emph{flat
  co-ordinates for $Y$}.  Note that the exponentiated flat
co-ordinates $\exp(g^k(y))$ are single-valued.

The $J$-functions satisfy differential equations which determine the
quantum products.

\begin{pro}
  \label{pro:QDE} \ 
  \begin{enumerate}
  \item 
    \[
    z \parfrac{}{\tau^i} z \parfrac{}{\tau^j} J_{\cX}(\tau,z) =
    \sum_{k=0}^{n-1} \(\delta_i \underset{\rm big}{\star}\)_j^{\phantom{j}k} 
    z \parfrac{}{\tau^k} J_{\cX}(\tau,z) 
    \]
    where $ \(\delta_i \underset{\rm big}{\star}\)_j^{\phantom{j}k} $ are the matrix
    entries of the product \eqref{eq:QCX}.
  \item 
    \[
    z \parfrac{}{t^i} z \parfrac{}{t^j} J_{Y}(t,z) =
    \sum_{k=0}^{n-1} \(\gamma_i \underset{\rm big}{\star}\)_j^{\phantom{j}k} 
    z \parfrac{}{t^k} J_{Y}(t,z) 
    \]
    where $ \(\gamma_i \underset{\rm big}{\star}\)_j^{\phantom{j}k} $ are the matrix
    entries of the product \eqref{eq:QCY}.
  \end{enumerate}
\end{pro}

\begin{proof}
  Part (2) is well-known: it follows from the Topological Recursion
  Relations (\emph{cf.}
  \citelist{\cite{Pandharipande}*{Proposition~2}
    \cite{Cox--Katz}*{Chapter~10}}).  The proof of (1) is essentially
  identical, but uses the Topological Recursion Relations for
  orbifolds \cite[Section~2.5.5]{Tseng} instead of the Topological
  Recursion Relations for varieties.  Details here are in
  \cite{CCIT:An}.
\end{proof}

\subsection*{From PF to QC}

Propositions \ref{pro:mirror} and \ref{pro:QDE} together show that we
can determine the quantum products \eqref{eq:QCX} and \eqref{eq:QCY}
by looking at the differential equations satisfied by $I_\cX$ and
$I_Y$ \emph{in flat co-ordinates}:
\begin{align}
  \label{eq:QDEIX}
  z \parfrac{}{\tau^i} z \parfrac{}{\tau^j} I_{\cX}(x(\tau),z) &=
  \sum_{k=0}^{n-1} \(\delta_i \underset{\rm big}{\star}\)_j^{\phantom{j}k} 
  z \parfrac{}{\tau^k} I_{\cX}(x(\tau),z) \\
  \label{eq:QDEIY}
  z \parfrac{}{t^i} z \parfrac{}{t^j} I_{Y}(y(t),z) &=
  \sum_{k=0}^{n-1} \(\gamma_i \underset{\rm big}{\star}\)_j^{\phantom{j}k} 
  z \parfrac{}{t^k} I_{Y}(y(t),z) 
\end{align}
A more invariant way to say this is as follows.  Let $\lambda_1$,
$\lambda_2$ be fixed complex numbers.  If we associate to a vector
field $v = \sum v_k(y) \parfrac{}{y_k}$ on $\cM_B$ the differential
operator $\sum z v_k(y) \parfrac{}{y_k}$ then the systems of
differential equations \eqref{eq:PFX}, \eqref{eq:PFY} define a
$\cD$-module on $\cM_B$.  The characteristic variety $\CharV$ of this
$\cD$-module is a subscheme of $T^\star \cM_B$, and we can read off
the quantum products from the algebra of functions $\cO_\CharV$.
Indeed, choosing flat co-ordinates on a neighbourhood $U$ of the
large-radius limit point for $\cX$ in $\cM_B$ identifies $\cO_U$ with
analytic functions in $\tau^0,\ldots,\tau^{n-1}$ and identifies the algebra
of fiberwise-polynomial functions on $T^\star U$ with
$\cO_U[\xi_0,\ldots,\xi_{n-1}]$; here $\xi_k$ is the fiberwise-linear
function on $T^\star U$ given by $\parfrac{}{\tau^k}$.  The ideal
defining the characteristic variety $\CharV$ is generated by elements
\[
P(\tau^0,\ldots,\tau^{n-1},\xi_0,\ldots,\xi_{n-1},0)
\]
where $P(\tau^0,\ldots,\tau^{n-1},\xi_0,\ldots,\xi_{n-1},z)$ runs over the
set of fiberwise-polynomial functions on $T^\star U$ which depend
polynomially on $z$ and satisfy
\[
P\(\tau^0,\ldots,\tau^{n-1},z\parfrac{}{\tau^0},\ldots,z\parfrac{}{\tau^{n-1}},z\)
I_\cX(u,z) = 0.
\]
Equation \eqref{eq:QDEIX} implies that
\[
\left.\cO_\CharV \right|_U = \cO_U[\xi_0,\ldots,\xi_{n-1}]/\mathfrak{I}
\]
where the ideal $\mathfrak{I}$ is generated by
\begin{align*}
\xi_i \xi_j -
\sum_{k=0}^{n-1} \(\delta_i \underset{\rm big}{\star}\)_j^{\phantom{j}k} \xi_k && 0 \leq i,j<n.
\end{align*}
In other words, the quantum cohomology algebra \eqref{eq:QCX} of $\cX$
is the algebra of functions $\left.\cO_\CharV \right|_U$ on the
characteristic variety $\CharV$, \emph{written in flat co-ordinates on
  $U$}.  

Similarly, choosing flat co-ordinates on a neighbourhood $V$ of the
large-radius limit point for $Y$ in $\cM_B$ identifies $\cO_V$ with
analytic functions in $t^0,\ldots,t^{n-1}$, and identifies the algebra
of fiberwise-polynomial functions on $T^\star V$ with
$\cO_V[\eta_0,\ldots,\eta_{n-1}]$ where $\eta_k$ is the
fiberwise-linear function on $T^\star V$ given by $\parfrac{}{t^k}$.
Equation \eqref{eq:QDEIY} implies that
\[
\left.\cO_\CharV \right|_V = \cO_V[\eta_0,\ldots,\eta_{n-1}]/\mathfrak{J}
\]
where the ideal $\mathfrak{J}$ is generated by
\begin{align*}
\eta_i \eta_j -
\sum_{k=0}^{n-1} \(\gamma_i \underset{\rm big}{\star}\)_j^{\phantom{j}k} \eta_k && 0 \leq i,j<n,
\end{align*}
and so the quantum cohomology algebra \eqref{eq:QCY} of $Y$ is the
algebra of functions $\left.\cO_\CharV \right|_V$ on the
characteristic variety $\CharV$, \emph{written in flat co-ordinates on
  $V$}.

\label{sec:discussion}

The characteristic variety $\CharV$ is a global analytic object ---
$\cO_\CharV$ gives an analytic sheaf of $\cO_{\cM_B}$-algebras,
defined over all of $\cM_B$ --- so to show that the quantum cohomology
algebras of $\cX$ and of $Y$ are related by analytic continuation
followed by the change-of-variables
\[
t^i =
\begin{cases}
  \tau^0, & i=0 \\
  -{2 \pi \sqrt{-1} \over n} + \sum_{j=1}^{n-1} L^i_{\phantom{i} j} \tau^j,  & i>0
\end{cases}
\]
we just need to show that the flat co-ordinates for $\cX$ and for $Y$
are related by analytic continuation followed by the
change-of-variables
\[
g^i(y) = 
\begin{cases}
  f^0(x), & i=0 \\
  -{2 \pi \sqrt{-1} \over n} + \sum_{j=1}^{n-1} L^i_{\phantom{i} j} f^j(x),  & i>0.
\end{cases}
\]
We will do this by showing that $g^i(y)$ and $f^j(x)$ satisfy the same
system of differential equations.

\subsection*{The GKZ System Associated to $Y$}

The \emph{GKZ system} associated to $Y$ is the system of differential
equations
\begin{multline}
  \label{eq:GKZ}
  \(
  \prod_{j:D'_j(d) >0} \prod_{m=0}^{D'_j(d) - 1} \(\gimel_j - m\) 
  \) f \\
  = y_1^{d_1} \cdots y_{n-1}^{d_{n-1}} \(
  \prod_{j:D'_j(d) <0} \prod_{m=0}^{-D'_j(d) - 1} \(\gimel_j - m\) 
  \) f
\end{multline}
for every $d = d_1 \beta_1 + \cdots + d_{n-1} \beta_{n-1}$, where
\[
\gimel_j =
\begin{cases}
  y_1 \parfrac{}{y_1} & j=0 \\
  -2 y_1 \parfrac{}{y_1} + y_2 \parfrac{}{y_2} & j=1 \\
  y_{j-1} \parfrac{}{y_{j-1}} - 2 y_j \parfrac{}{y_j} + y_{j+1} \parfrac{}{y_{j+1}} & 1<j<n-1 \\
  y_{n-2} \parfrac{}{y_{n-2}} - 2 y_{n-1} \parfrac{}{y_{n-1}} & j=n-1 \\
  y_{n-1} \parfrac{}{y_{n-1}} & j=n.
\end{cases}
\]

\begin{pro}
  Both
  \begin{itemize}
  \item[(a)] $f^1(x), \ldots, f^{n-1}(x)$ plus the constant function; and
  \item[(b)]   $g^1(y),\ldots,g^{n-1}(y)$ plus the constant function
  \end{itemize}
  form bases of solutions to the GKZ system \eqref{eq:GKZ}.
\end{pro}

\begin{proof}
  The sets (a) and (b) are evidently each linearly independent.  The
  constant function evidently satisfies \eqref{eq:GKZ}.  The flat
  co-ordinates $f^1(x), \ldots, f^{n-1}(x)$ and
  $g^1(y),\ldots,g^{n-1}(y)$ are independent of $\lambda_1$,
  $\lambda_2$, $x_0$, and $y_0$, so they can be extracted from the
  $z^0$ terms of $I_\cX$ and $I_Y$ after setting $\lambda_1 =
  \lambda_2 = x_0 = y_0 = 0$.  But $\left. I_\cX \right|_{\lambda_1 =
    \lambda_2 = x_0 = 0}$ and $\left. I_Y \right|_{\lambda_1 =
    \lambda_2 = y_0 = 0}$ satisfy the systems of differential
  equations \eqref{eq:PFXa}, \eqref{eq:PFYa} with $\lambda_1$ and
  $\lambda_2$ set to zero, and once $\lambda_1$ and $\lambda_2$ are
  set to zero the $z$-dependence in these differential equations
  cancels.  The resulting system of differential equations in each
  case is \eqref{eq:GKZ}.
\end{proof}

Any analytic continuation $\tilde{g}_i(y)$ of $g^i(y)$ to a
neighbourhood of the large-radius limit point for $\cX$ still
satisfies \eqref{eq:GKZ}, so
\[
\tilde{g}_i(y) = \sum_{j=1}^{n-1} L^i_{\phantom{i} j} f^j(x) + m_i
\]
for some constants $L^i_{\phantom{i} j}$ and $m_i$.  Thus any analytic
continuation of $g^i(y)$ is an affine-linear combination of the flat
co-ordinates $f^1(x),\ldots,f^{n-1}(x)$.  It remains to choose a
specific analytic continuation and determine the corresponding
constants $L^i_{\phantom{i} j}$ and $m_i$.  Before we do so, we
prove Proposition~\ref{pro:inversemirrorA_n}.  This follows
immediately from:

\begin{pro} \label{pro:flatcoordsarebasis}
  Let $\zeta = \exp\({\pi \sqrt{-1} \over n}\)$.  Let
  $\kappa_0(x),\ldots,\kappa_{n-1}(x)$ be the roots of the polynomial
  \[
  W_\cX(\kappa) = \kappa^n + x_{n-1} \kappa^{n-1} + x_{n-2}
  \kappa^{n-2} + \cdots + x_{1} \kappa + 1,
  \]
  where the roots are numbered such that as $x \to 0$, $\kappa_i(x)
  \to \zeta^{2i+1}$, $0 \leq i < n$.  Then another basis of solutions
  to the GKZ system \eqref{eq:GKZ} is given by
  \begin{align} \label{eq:logdifferencekappai}
    \log \kappa_i(x) - \log \kappa_{i-1}(x), &&
    1 \leq i < n
  \end{align}
  together with the constant function.  Furthermore, 
  \begin{equation}
    \label{eq:roottof}
    \log \kappa_i(x) = {(2i+1) \pi \sqrt{-1} \over n} + {1 \over n}
    \sum_{k=1}^{n-1} \zeta^{(2i+1) k} f^k(x).
  \end{equation}
\end{pro}

\begin{proof}
  Most of this is classical: see for example
  \citelist{\cite{Stienstra}*{Theorem~1}\cite{Mayr}}.  Here we follow
  \cite[Section~6]{Hori--Vafa}.  Let $(\kappa,\nu,u,v)$ be
  co-ordinates on $\Cstar \times \CC^3$. By using the Morse--Bott
  function $(\kappa,\nu,u,v)\mapsto \Re(\nu(W_\cX(\kappa)+uv))$ on
  $\Cstar\times \CC^3$, we can construct a Morse cycle $\Gamma\subset
  \Cstar\times \CC^3$ as the union of downward gradient flowlines for
  $\Re(\nu(W_\cX(\kappa)+uv))$ from a compact 2-cycle $\Gamma'$ in the
  critical set $Z:=\{\nu=0, W_\cX(\kappa)+uv=0\}\subset \Cstar\times
  \CC^3$.  The integral
  \[
  \gamma(x) = \int_\Gamma \exp\( {\nu (W_\cX(\kappa) + uv)} \)
  {d\kappa \over \kappa} {d\nu} \, du \, dv
  \]
  satisfies the differential equations \eqref{eq:GKZ}.  But the
  integrand here is equal to
  \[
  d\( {\exp\({\nu  \(W_\cX(\kappa) + u v\)}\) \over W_\cX(\kappa) + uv} {d\kappa \over \kappa} du \, dv\)
  \]
  outside $Z$, and so
  \begin{align*}
    \gamma(x) &= \int_{\Gamma' \subset Z} {{d\kappa \over \kappa} du \, dv \over
      d\(W_\cX(\kappa) + u v\)} \\
    &= \int_{\Gamma' \subset X} {d\kappa \over \kappa}{du \over u}.
  \end{align*}
  Now we take $\Gamma'\subset X$ to be a vanishing cycle for the
  function $W_\cX + uv$.  Integrating out $du/u$ gives
  \begin{equation}
    \label{eq:differencelog}
    \gamma(x) = \int_{C\subset \Cstar} \frac{d\kappa}{\kappa} = \log \kappa_m - \log \kappa_l,    
  \end{equation}
  where $\kappa_l,\kappa_m$ are roots of $W_\cX(\kappa)=0$ and $C$ is
  a path from $\kappa_l$ to $\kappa_m$.  By choosing an appropriate
  basis of vanishing cycles $\Gamma'$, we find the $n-1$ linearly
  independent solutions \eqref{eq:logdifferencekappai}; these,
  together with the constant function, form a basis of solutions to
  the GKZ system.
  
  It remains to prove \eqref{eq:roottof}.  We have
  \[
  \kappa_i(x) = \zeta^{2i+1} +O\(x_j\)  
  \]
  and 
  \begin{align*}
    x_{k} = (-1)^{n-k} e_{n-k}(\kappa_0,\dots,\kappa_{n-1}),
    &&  1 \leq k < n.
  \end{align*}
  Since the constant function and \eqref{eq:differencelog} are
  solutions to the GKZ system \eqref{eq:GKZ} and $\kappa_0\cdots
  \kappa_{n-1} = (-1)^n$, each $\log \kappa_i$ is also a solution to
  \eqref{eq:GKZ}.  Thus, by Proposition~\ref{pro:flatcoordsarebasis},
  $\log \kappa_i$ must be a linear combination of
  $f^1(x),\dots,f^{n-1}(x)$ and a constant:
  \[
  \log \kappa_i = \sum_{r=1}^{n-1} c_{ir} f^r(x) + \text{const}.    
  \]
  Since $\kappa_i=\zeta^{2i+1}$ at $x=0$, we have 
  \begin{equation}
    \label{eq:logx}
    \log \(\frac{\kappa_i}{\zeta^{2i+1}}\) = \sum_{r=1}^{n-1} c_{ir} f^r(x).
  \end{equation}
  On the other hand, differentiating $W_\cX(\kappa_i)=0$ gives
  \[
  (\kappa_i)^k + n (\kappa_i)^{n-1} \parfrac{\kappa_i}{x_k} + O\(x_j\) = 0
  \]
  and so
  \[
  \kappa_i(x) = \zeta^{2i+1}  - {1 \over n} \sum_{k=1}^{n-1}
  \zeta^{(2i+1)(k-n+1)} x_k  +O\Big(\(x_j\)^2\Big).  
  \]
  Substituting this into \eqref{eq:logx} and using the fact that 
  \[
  f^r(x) = x_r + O\Big(\(x_j\)^2\Big)
  \]
  yields \eqref{eq:roottof}.
\end{proof}

We observed above Proposition~\ref{pro:flatcoordsarebasis} that any
analytic continuation of $g^i(y)$ is an affine-linear combination of
the flat co-ordinates $f^1(x),\ldots,f^{n-1}(x)$.

\begin{pro} \label{pro:ac}
  There exists a path from the large-radius limit point for $Y$ to the
  large-radius limit point for $\cX$ such that the analytic
  continuation of the flat co-ordinates $g^i(y)$, $1 \leq i<n$, along
  that path satisfy
  \[
  g^i(y) = - {2 \pi \sqrt{-1} \over n} + {1 \over n} \sum_{k=1}^{n-1}
  \zeta^{2 k i} \( \zeta^{-k }- \zeta^{k}\) f^k(x),
  \]
  where $\zeta = \exp\({\pi \sqrt{-1} \over n}\)$.
\end{pro}

\begin{proof}
  Consider the polynomial
  \begin{multline*}
    W_Y(\mu) = \mu^n + \mu^{n-1} + y_1 \mu^{n-2} +
    y_1^2 y_2 \mu^{n-3} + y_1^3 y_2^2 y_3 \mu^{n-4} + \cdots
    \\
    + y_1^{n-1} y_2^{n-2} \cdots y_{n-2}^2 y_{n-1}
  \end{multline*}
  and number its roots $\mu_i(y)$, $0 \leq i<n$ such that as $y\to
  0$ 
  \begin{align*}
    \mu_0(y) & \to -1 \\
    \mu_1(y) & \sim -y_1 \\
    \mu_2(y) & \sim -y_1 y_2 \\
    & \vdots \\
    \mu_{n-1}(y) & \sim -y_1 y_2 \cdots y_{n-1}.\\
  \end{align*}
  We have $W_\cX(\kappa) = 0$ if and only if $W_Y(1/(x_1 \kappa)) = 0$,
  where $x_i$ and $y_j$ are related by \eqref{eq:xtoy}, so still
  another basis of solutions to the GKZ system \eqref{eq:GKZ} is
  \begin{align*}
    \log \mu_i(y) - \log \mu_{i-1}(y),&& 1 \leq i < n,
  \end{align*}
  together with the constant function.  The solution $g^i(y)$ is
  singled out by its behaviour $g^i(y) = \log y_i +
  O(y_1,\ldots,y_{n-1})$ as $y \to 0$, so
  \[
  g^i(y) = \log \mu_i(y) - \log \mu_{i-1}(y).
  \]
  
  Along any path from the large-radius limit point for $Y$ to the
  large-radius limit point for $\cX$, the root $\mu_i(y)$ of $W_Y$
  analytically continues to the root $1/(x_1
  \kappa_{\sigma(i)}(x))$ of $W_{\cX}$, for some permutation
  $\sigma$ of $\{0,1,\ldots,n-1\}$.  The group of monodromies around
  the discriminant locus of $W_\cX$ acts $n$-transitively on the set
  of roots of $W_\cX$, so we can choose a path such that $\sigma$ is
  the identity permutation.  Along this path, $\log \mu_i(y) - \log
  \mu_{i-1}(y)$ analytically continues to $\log \kappa_{i-1}(x) -
  \log \kappa_{i}(x)$, $1 \leq i < n$.  Applying equation
  \eqref{eq:roottof} yields
  \[
    g^i(y) = -{2 \pi \sqrt{-1} \over n} + {1 \over n} \sum_{k=1}^{n-1}
  \zeta^{2 k i} \( \zeta^{-k} - \zeta^{k}\) f^k(x).
  \]
\end{proof}

\begin{rem*}
  For an explicit path satisfying the conditions in
  Proposition~\ref{pro:ac}, we can concatenate two paths defined
  as follows.  The first runs from $(y_0,y_1, \ldots, y_{n-1}) =
  (0,0,\ldots,0,0)$ to $(y_0,y_1, \ldots, y_{n-1}) = (0,1,1,\ldots,1,1)$
  and is given by $y_0 = 0$ and
  \begin{align*}
    W_Y(\mu) = \Big(\mu - \(-1 - \epsilon \rho^2 - \epsilon^2 \rho^3 -
    \ldots - \epsilon^{n-1} \rho^n\)\Big) \prod_{k=1}^{n-1} \(\mu -
    \epsilon^k \rho^{k+1}\), && 0 \leq \epsilon \leq 1,
  \end{align*}
  where $\rho = \exp\({2 \pi \sqrt{-1} \over n+1}\)$.  The second runs
  from $(x_0,x_1, \ldots, x_{n-1}) = (0,1,1,\ldots,1,1)$ to $(x_0,x_1,
  \ldots, x_{n-1}) = (0,0,\ldots,0,0)$, and is given by $x_0 = 0$ and
  \begin{align*}
    W_\cX(\kappa) = \prod_{k=0}^{n-1} \( \kappa - 
    \exp \( \pi \sqrt{-1} \left[ {2k+1 \over n} \epsilon' + {2(n-k) \over
        n+1} (1 - \epsilon') \right] \) \), &&
    0 \leq \epsilon' \leq 1.
  \end{align*}
  Note that the points $(y_0,y_1, \ldots, y_{n-1}) = (0,1,1,\ldots,1,1)$
  and $(x_0,x_1, \ldots, x_{n-1}) = (0,1,1,\ldots,1,1)$ coincide.
\end{rem*}

\subsection*{The Proof of Theorem \ref{thm:maintheorem}}
 
Combining Proposition~\ref{pro:ac} 
with the discussion above equation
(\ref{eq:GKZ}) shows that 
the quantum cohomology algebras of $\cX$ and
$Y$ coincide after analytic continuation 
along the path specified in
Proposition~\ref{pro:ac} followed by the affine-linear
change-of-variables
\begin{align*}
  t^i &=
  \begin{cases}
    \tau^0, & i=0 \\
    -{2 \pi \sqrt{-1} \over n} + \sum_{j=1}^{n-1} L^i_{\phantom{i} j} \tau^j,  & i>0,
  \end{cases}
 &  L^i_{\phantom{i} j} & = {\zeta^{2 i j} \( \zeta^{-j} - \zeta^{j}\) \over n},
\end{align*}
and the linear isomorphism 
\begin{equation*}
  \begin{aligned}
    L:\HorblocTX & \to \HlocTY \\
    \delta_0 & \mapsto \gamma_0,\\
    \delta_j & \mapsto \sum_{i=1}^{n-1} L^i_{\phantom{i} j} \gamma_i, && \qquad 1 \leq j < n.
  \end{aligned}
\end{equation*}
To see that $L$ preserves the Poincar\'e pairings, first observe that
the bases
\begin{align*}
n \lambda_1 \lambda_2, \gamma_1,\gamma_2,\ldots,\gamma_{n-1} 
&&\text{and}&&
1, \omega_1,\omega_2,\ldots,\omega_{n-1} 
\end{align*}
for $\HlocTY$ are dual with respect to the Poincar\'e pairing on
$\HlocTY$.  Let $L^\dagger$ denote the adjoint to $L$ with respect to
the Poincar\'e pairing $(\cdot,\cdot)_Y$ and the orbifold Poincar\'e
pairing $(\cdot,\cdot)_\cX$.  It suffices to show that $(L^\dagger
\gamma,L^\dagger \gamma')_\cX = (\gamma,\gamma')_Y$ for all $\gamma,
\gamma' \in \HlocTY$.  For $1 \leq i<n$, we have $(L^\dagger
\omega_i,\delta_k)_\cX = (\omega_i,L \delta_k)_Y = L^i_{\phantom{i} k}$, and so 
\begin{align*}
L^\dagger \omega_i = n \sum_{k=1}^{n-1} L^i_{\phantom{i} k} \delta_{n-k}, && 1 \leq
i<n.
\end{align*}
Also $L^\dagger 1 = \delta_0$.  Straightforward calculation now gives
$(L^\dagger 1, L^\dagger 1)_\cX = (n \lambda_1 \lambda_2)^{-1}$,
$(L^\dagger 1, L^\dagger \omega_i)_\cX = 0$ for $1 \leq i<n$, and
\begin{align*}
  (L^\dagger \omega_i, L^\dagger \omega_j)_\cX =
  \begin{cases}
    0 & \text{if $|i-j|>1$} \\
    1 & \text{if $|i-j|=1$}\\
    -2 & \text{if $i=j$}
  \end{cases}
  && \text{for $1 \leq i,j<n$.}
\end{align*}
As the class $\omega_j$ is the $T$-equivariant Poincar\'e-dual to the
$j$th exceptional divisor we see that $L^\dagger$, and hence $L$, is
pairing-preserving.  This completes the proof of
Theorem~\ref{thm:maintheorem}.  \hfill $\Box$

\begin{rem*}
  A more conceptual explanation of this result is as follows.  One can
  construct a Frobenius manifold from a \emph{variation of
    semi-infinite Hodge structure} \cite{Barannikov} (henceforth \VHS)
  together with a choice of \emph{opposite subspace}\footnote{Mirror
    symmetry often associates to the quantum cohomology of some target
    space a ``mirror family'' of manifolds.  In this case one can
    think of the \VHS as an analog of the usual variation of Hodge
    structure on the mirror family, and the opposite subspace as an
    analog of the weight filtration.}.  We have argued elsewhere that
  in certain toric examples one can construct the Frobenius manifold
  which is the ``mirror partner'' to the quantum cohomology of $Y$
  from a \VHS parameterized by the $B$-model moduli space of $Y$,
  together with a distinguished opposite subspace associated to the
  large-radius limit point for $Y$ \cite{CCIT:wallcrossings1}.  (The
  Frobenius manifold mirror to the quantum cohomology of a toric
  orbifold $\cX$ birational to $Y$ is given by the \emph{same} \VHS
  but the opposite subspace corresponding to the large-radius limit
  point for $\cX$.)  One can apply this construction here to get a
  \VHS parametrized by $\cM_B$.  This \VHS has the special property
  that the opposite subspace at the large-radius limit point for $Y$
  agrees with the opposite subspace at the large-radius limit point
  for\footnote{In fact each maximal cone of the secondary fan gives
    rise to a toric partial resolution $Y'$ of $\cX$, a large-radius
    limit point for $Y'$, and an opposite subspace corresponding to
    this large-radius limit point.  All these opposite subspaces
    agree: we really do get a Frobenius structure defined over all of
    $\cM_B$.}
  $\cX$.  In general the
  difference between the opposite subspaces at different large radius
  limit points will be measured by an element of Givental's linear
  symplectic group, but in this case the corresponding group element
  maps the opposite subspaces isomorphically to each other.  This
  means that we get a Frobenius manifold \emph{over the whole
    (non-linear) space $\cM_B$}.  One can construct flat co-ordinates
  in a neighbourhood of any point of $\cM_B$, and the transition
  functions between such flat co-ordinate patches, such as
  \[
  g^i(y) = \sum_j L^i_{\phantom{i} j} f^j(x) + \log c_i,
  \]
  are necessarily affine-linear and (Poincar\'e) metric-preserving.
\end{rem*}

\begin{rem*}
  It is clear from the proof of Proposition~\ref{pro:ac} that changing
  the path along which analytic continuation is taken will result in a
  corresponding change in the statements of
  Theorem~\ref{thm:maintheorem} and its Corollary.  Hence the
  co-ordinate change in Theorem~\ref{thm:maintheorem} is not unique.
  This ambiguity can be understood as an automorphism of quantum
  cohomology.  The orbifold fundamental group
  \[
  G := \pi_1^{\text{\rm orb}}\(\cM_B\setminus\{\text{discriminant locus of
    $W_\cX$}\}\)
  \]
  acts simply-transitively on the set of homotopy types of paths from
  the large-radius limit point for $Y$ to that for $\cX$, and in
  particular acts transitively (but ineffectively) on the set of all
  possible co-ordinate changes obtained by analytic continuation.
  This deserves further study: we just note here the intriguing fact
  that $G$ is isomorphic to $\widetilde{A}_{n-1}\rtimes \ZZ_n$, which
  also appears as a subgroup (generated by spherical twists and line
  bundles) of the group of autoequivalences of $D_E^b(Y)$ \cite{IU,
    Bridgeland, IUU}.  Here $\widetilde{A}_{n-1}$ is the affine braid
  group and $D_E^b(Y)$ is the bounded derived category of coherent
  sheaves on $Y$ supported on the exceptional set $E$.
\end{rem*}

\section{Givental's Lagrangian Cone as a Formal Scheme}
\label{sec:rigorouscone}

In the main body of the text, following Givental, we encode genus-zero
Gromov--Witten invariants via the formal germ $\cL$ ($= \cLs\text{ or
}\cL^\untw$) of an infinite-dimensional submanifold in a symplectic
vector space $\cH$.  This submanifold-germ $\cL$ has very special
geometric properties, several of which play an essential role in the
proof of Theorem~\ref{thm:generaltwist}.  These
geometric properties were discovered in \cite{Coates--Givental:QRRLS,
  Givental:symplectic}.  The properties of $\cL$ described in
\emph{loc. cit.} should, as stated there, ``be interpreted in the
sense of formal geometry''.  It is clear what this means, but an
appropriate framework of formal geometry seems to be missing from the
literature.  In this Appendix we remedy this: we define the formal
germ $\cL$ and the formal germ of $\cH$ as non-Noetherian formal
schemes.  We also establish, within our framework, the geometric
properties of $\cL$ which we use above.  A more detailed account of
the geometry described here will be given in a future paper
\cite{CCIT:cone}.

This section is unavoidably technical, and we urge the reader
unfamiliar with Givental's formalism to begin by reading the much more
accessible account \cite{Givental:symplectic}.  Also, we should
emphasize that any originality here is in the framework that we use
rather than in the content of our arguments.  The geometric properties
of $\cL$ which we discuss have already been established by Givental
in \cite{Givental:symplectic}, and our arguments share a core of ideas
with his.

\subsection*{Topological rings} 
In the main body of the text, we considered topological rings equipped
with a non-negative additive valuation.  Here we will deal with a
broader class of topological rings.  Let $R$ be a commutative
topological ring with a unit.  A topology on $R$ is said to be
\emph{linear} if a fundamental neighborhood system of zero in $R$ is
given by a descending chain of ideals $I_0 \supset I_1 \supset I_2
\supset \cdots$ in $R$.  Conversely, any descending chain
$\{I_n\}_{n\ge 0}$ of ideals defines a unique linear topology on $R$.
In this case we say that the topology on $R$ \emph{is defined by
  ideals} $\{I_n\}_{n\ge 0}$.  For example, a non-negative additive
valuation $v \colon R \setminus \{0\} \to \R_{\ge 0}$ on an integral
domain $R$ gives a linear topology on $R$ defined by the ideals $I_n
=\{ r\in R \;;\; v(r)\ge n\}$.  Throughout Appendix
\ref{sec:rigorouscone}, \emph{a topology on a ring is assumed to be
  linear, complete and Hausdorff}\footnote{ In the literature on
  formal schemes, an additional admissibility condition is imposed on
  $R$.  For example, McQuillan \cite{McQuillan} considered the
  condition that $R^{\rm nilp}$ is open, \emph{i.e.} contains some
  $I_n$.  Although our examples always satisfy this condition, we do
  not need this in the explanation below.}.

Define the space of \emph{convergent Laurent series} in $z$ to be
\[
R\{z,z^{-1}\} := \left
\{\sum_{n\in \Z} r_n z^n \;:\; r_n\in R, \; 
r_n \to 0 \text{ as } |n| \to \infty
\right \}. 
\]
When the topology on $R$ is defined by $\{I_n\}_{n\ge 0}$, this space
is identified with the completion of $R[z,z^{-1}]$ with respect to the
topology defined by $\{I_n[z,z^{-1}]\}_{n\ge 0}$.  When the topology
on $R$ comes from an additive valuation, this definition coincides
with that given in Section~\ref{sec:GiventalFormalism}.  

Define the ideal of \emph{topologically nilpotent elements} to be
\[
R^{\rm nilp}:=
\left\{r\in R\;:\; \lim_{n\to \infty} r^n =0
\right\}. 
\]
When the topology on $R$ is defined by ideals $\{I_n\}_{n\ge 0}$, the
topology on $R[\epsilon]/\langle \epsilon^2\rangle$ is defined by
$\{I_n + I_n \epsilon\}_{n\ge 0}$, that on $R[\![t]\!]$ is defined by
$\{I_n [\![t]\!] + t^n R[\![t]\!]\}_{n \ge 0}$, and that on
$R\{z,z^{-1}\}$ is given by $\{I_n\{z,z^{-1}\}\}_{n\ge 0}$.  One can
check that if $R$ is complete and Hausdorff then so are
$R[\epsilon]/\langle \epsilon^2\rangle$, $R[\![t]\!]$, and
$R\{z,z^{-1}\}$.

\subsection*{Notation Conventions}

Recall that $\{\phi_\alpha : 1\le \alpha \le N\}$ and $\{\phi^\alpha:
1\le \alpha \le N\}$ are bases for $H_{\rm orb}^\bullet(\cX)$ such
that $(\phi_\alpha,\phi^\beta)_{\rm orb}^\tw = {\delta_\alpha}^\beta$.
We take $\phi_1$ to be the unit class $\fun$, and set $g_{\alpha\beta}
= (\phi_\alpha,\phi_\beta)_{\rm orb}^\tw$ and $g^{\alpha\beta} =
(\phi^\alpha,\phi^\beta)_{\rm orb}^\tw$.  Throughout this Appendix we
use Einstein's summation convention for Greek indices, summing
repeated Greek (but not Roman) indices over the range $1,2,\ldots,N$.

\subsection*{A formal germ at $-z$ of Givental's space $\cH$} 
Recall that $\Lambda[\![s_0,s_1,\dots]\!]$ was defined to be the
completion of $\C[\EffX][s_0,s_1,\dots]$ with respect to the
non-negative additive valuation $v$ given by
\[
v(Q^d) = \int_d \omega, \quad v(s_k) = k+1. 
\]
Here the $s_i$ are parameters of the universal invertible
multiplicative characteristic class $\bc$: see
(\ref{eq:universalclass}).  We write $\Lambda_\bs :=
\Lambda[\![s_0,s_1,\dots]\!]$.  Let $\cH$ be Givental's vector space
over $\Lambda_\bs$:
\[
\cH = H^\bullet_{\rm orb}(\cX;\C) \otimes_\C \Lambda_\bs\{z,z^{-1}\}. 
\]
We write a general element $f\in \cH$ in the form:  
\begin{align*}
  f = -z + \bt(z) + \bp(z), &&
  \bt(z) = \sum_{k\ge 0} t^\alpha_k 
  \phi_\alpha z^k, &&
  \bp(z) = \sum_{l \ge 0} p_{l\beta}
  \frac{\phi^\beta}{(-z)^{l+1}}. 
\end{align*}
We regard the coefficients $t_k^\alpha, p_{l\beta}\in \Lambda_\bs$ as
affine co-ordinate functions $t_k^\alpha,p_{l\beta}\colon \cH\to
\Lambda_\bs$.  We define the formal germ $(\cH,-z)$ at $-z$ to be the
affine formal scheme over $\Lambda_\bs$:
\begin{align*}
  (\cH,-z) := \Spf S, &&
  S:= \Lambda_\bs[ t_k^\alpha,
  p_{l\beta} : 
  1\le \alpha,\beta \le N; 0 \leq k,l<\infty ]\sphat,  
\end{align*}
where $\sphat$\ means the completion with respect to the valuation
\[
v(t_{k}^\alpha) = k+1, \quad 
v(p_{l\beta}) = l+1. 
\]
We regard the formal scheme $(\cH,-z) = \Spf(S)$ as a \emph{functor}
from the category of topological $\Lambda_\bs$-algebras and continuous
$\Lambda_\bs$-algebra homomorphisms to the category of sets:
\begin{align*}
  (\cH,-z) \colon \text{(Topological $\Lambda_\bs$-algebras)} &
  \longrightarrow \text{(Sets)} \\
  R & \longmapsto \Hom (S, R).
\end{align*}
We have 
\[
(\cH,-z)(R) \cong \left\{
-z + \sum_{n\in \Z} r_n^\alpha \phi_\alpha z^n\;:  
\; r_n^\alpha \in R^{\rm nilp}, \;   
\;  r_n^\alpha \to 0  
\text{ as } |n| \to \infty\right \}.  
\]
Using this identification, we always write an element of $(\cH,-z)(R)$
in the form $f= -z + \sum_{n\in \Z} r_n^\alpha \phi_\alpha z^n = -z +
\bt(z) + \bp(z)$ with $t_k^\alpha =r_k^\alpha$, $p_{l\beta} =
(-1)^{l+1} r_{-l-1}^{\nu}g_{\nu\beta}$.  The tangent functor
$T(\cH,-z)$ is defined to be
\[
T(\cH,-z)(R) := 
(\cH,-z)(R[\epsilon]/\langle \epsilon^2\rangle). 
\]
The tangent space $T_f(\cH,-z)(R)$ at $f\in (\cH,-z)(R)$ is defined to
be the preimage of $f$ under the natural map 
$(\cH,-z)(R[\epsilon]/\langle \epsilon^2\rangle) \to (\cH,-z)(R)$.  We
have, for $f\in (\cH,-z)(R)$,
\begin{align*}
T_f(\cH,-z)(R) &\cong 
\left\{ \sum_{n\in \Z} v_n^\alpha \phi_\alpha 
z^n \;:\; v_n^\alpha \in R, \ v_n^\alpha \to 0 \text{ as } 
|n| \to \infty \right\} \\ 
&\cong H^\bullet_{\rm orb}(\cX;\C)\otimes_\C R\{z,z^{-1}\}. 
\end{align*} 
The tangent space $T_f(\cH,-z)$ is equipped with the topology induced
from that on $R\{z,z^{-1}\}$.  This coincides with the ``pointwise
convergence topology'' on $\Hom(S, R[\epsilon]/\langle
\epsilon^2\rangle) =T(\cH,-z)(R)$, \emph{i.e.}  $\phi_n \in
\Hom(S,R[\epsilon]/ \langle \epsilon^2\rangle )$ converges to
$\phi_\infty$ if and only if $\phi_n(s)$ converges to
$\phi_{\infty}(s)$ in $R[\epsilon]/\langle \epsilon^2\rangle$ for
every $s\in S$.

\subsection*{The Formal Subscheme $\cLs$ of $(\cH,-z)$} 
We consider genus zero Gromov-Witten theory twisted by the universal
multiplicative characteristic class $\bc$ in
(\ref{eq:universalclass}).  Introduce the double correlator notation:
\begin{multline*}
  \corrr{\phi_{\alpha_1}\psi^{k_1}, 
    \dots,\phi_{\alpha_m} \psi^{k_m}}_{\bt}^{\tw}   
  :=\\ \sum_{d\in \EffX} \sum_{n\ge 0}  \frac{Q^d}{n!} 
  \corr{\phi_{\alpha_1}\psi^{k_1},\dots,\phi_{\alpha_m}\psi^{k_m}, 
    \bt(\psi),\dots,\bt(\psi)}_{0,m+n,d}^{\cX, \tw},   
\end{multline*}
where the summand is defined to be zero in the unstable range $d=0$,
$m+n<3$.  These correlators are elements of the completion of
$\Lambda_\bs[t_k^\alpha: 1\le \alpha\le N; 0 \leq k < \infty]$.

\begin{exait*}
  The correlator $\corrr{\,\vphantom{\psi^k}}_{\bt}$ with no insertion
  is the genus-zero descendant potential.
\end{exait*}

Define elements $E_{j\alpha}\in S$ and the ideal $\mathfrak{I} \subset
S$ by
\begin{align*}
  E_{j\alpha} &= p_{j\alpha} - 
  \corrr{\phi_\alpha\psi^j}^{\tw}_{\bt}, &
  \mathfrak{I} &= \big\langle E_{j\alpha}:
  1\le \alpha\le N;0 \leq j < \infty \big\rangle.
\end{align*}
The formal subscheme $\cLs$ of $(\cH,-z)$ 
is defined to be 
\begin{equation} 
\label{eq:def_rigorouscone}
\cLs := \Spf\Big(S \big /\, \ov{\mathfrak{I}} \Big),   
\end{equation} 
where $\ov{\mathfrak{I}}$ is the closure of the ideal $\mathfrak{I}$.
Considering the virtual dimension of the moduli space of stable maps
shows that $E_{j\alpha}$ converges to zero in $S$ as $j\to \infty$.
Therefore the set $\cLs(R)$ of $R$-valued points of $\cLs$ is given by
the graph of the functions $\bt \mapsto \corrr{\phi_\alpha
  \psi^j}^{\tw}_{\bt}$:
\[
\cLs(R) = \left 
\{ -z + \bt(z) + \bp(z) 
\in (\cH,-z)(R)\;: \;  
p_{k\alpha} = \corrr{\phi_\alpha \psi^k}_{\bt}^\tw 
\; \forall k, \alpha \right \}. 
\]

\begin{exadefit*}  
  Let $\Lambda_\bs[\![\tau]\!]  :=
  \Lambda_\bs[\![\tau^1,\dots,\tau^N]\!]$.  The twisted $J$-function
  (\ref{eq:twistedJfun}) is a $\Lambda_\bs[\![\tau]\!]$-valued point
  on $\cLs$:
  \[
  J^\tw(\tau,-z) = -z + \tau + \sum_{k\ge 0} 
  \corrr{\phi_\alpha \psi^k}_{\bt = \tau}^\tw 
  \frac{\phi^\alpha}{(-z)^{k+1}} 
  \in \cL^\tw(\Lambda[\![\tau]\!])
  \]
  characterized by the condition (\ref{eq:Jcharact}). 
\end{exadefit*}

The tangent functor $T\cLs$, the tangent spaces $T_f\cLs(R)$, and the
topologies on them are defined as above.  Explicitly, the tangent
space $T_f\cLs(R)$ at $f = -z + \bt(z) + \bp(z)\in \cLs(R)$ is given
by the set of points $\sum_k \dot{t}_k^\alpha \phi_\alpha z^k + \sum_k
\dot{p}_{k\alpha} \frac{\phi^\alpha}{(-z)^{k+1}}$ in $H^\bullet_{\rm
  orb}(\cX)\otimes R\{z,z^{-1}\}$ satisfying
\begin{equation}
\label{eq:tangentvector}
\dot{p}_{k\alpha} = \sum_{l,\beta} \dot{t}_l^\beta 
\corrr{\phi_\beta \psi^l, \phi_\alpha\psi^k}_{\bt}^\tw. 
\end{equation} 
It is easy to check that $T_f\cLs(R)$ is a closed subspace of
$T_f(\cH,-z)= H^\bullet_{\rm orb}(\cX)\otimes R\{z,z^{-1}\}$.

The following elementary fact will be useful:
\begin{lem} 
\label{lem:diff_arc} 
If $I(t)\in \cLs(R[\![t]\!])$ then the derivative ${d I
  \over d t}(t)$ lies in $T_{I(t)} \cLs(R[\![t]\!])$.
\end{lem} 
\begin{proof} 
Observe that there exists an automorphism of 
$R[\![t]\!][\epsilon]/\langle \epsilon^2 \rangle$ 
which sends $t$ to $t+\epsilon$. 
Therefore, $I(t+\epsilon) = I(t) + \epsilon {dI \over dt}(t)$ 
belongs to $\cLs(R[\![t]\!][\epsilon]/\langle \epsilon^2\rangle)$. 
Thus ${dI \over dt}(t) \in T_{I(t)} \cLs(R[\![t]\!])$. 
\end{proof} 

\subsection*{Special geometric properties of $\cL=\cL^\untw$}

The special case $\bc=1$ (\emph{i.e.} $s_0=s_1=s_2=\cdots=0$) gives a
formal scheme $\cL^\untw$, defined over $\Lambda$, which corresponds
to untwisted Gromov-Witten theory.  We now verify some basic geometric
properties of $\cL=\cL^\untw$ which are used in the main body of the
text.  In the rest of the section we omit the superscript ``un'',
writing $\cL$ for $\cL^\untw$, $\corrr{\cdots\vphantom{\psi^k}}_\bt$
for $\corrr{\cdots\vphantom{\psi^k}}^\untw_\bt$, etc.  Also we take $R$
to be a complete, Hausdorff, linearly-topologized $\Lambda$-algebra.

The geometric properties of $\cL$ follow from the three universal
relations in genus-zero orbifold Gromov-Witten theory \cite[Section
3.1]{Tseng}:
\begin{align}
&\corrr{\psi \vphantom{\psi^k}}_{\bt}
= \sum_{k\ge 0} t_k^\alpha \corrr{\phi_\alpha \psi^k}_\bt
-2 \corrr{\vphantom{\psi^k}}_\bt,  \tag{DE} \\ 
&\corrr{\fun \vphantom{\psi^k}}_\bt  
= \frac{1}{2}(t_0,t_0)_{\rm orb}  
+ \sum_{k\ge 0} t_{k+1}^\alpha \corrr{\phi_\alpha \psi^k}_{\bt},
\tag{SE} \\ 
&\corrr{\phi_\alpha \psi^{k+1},
\phi_\beta\psi^l, \phi_\gamma \psi^{m}}_{\bt}
= \corrr{\phi_\alpha \psi^k, \phi_\nu}_{\bt} 
\corrr{\phi^\nu, \phi_\beta \psi^l, 
\phi_\gamma \psi^m}_{\bt}.  \tag{TRR}
\end{align}
These are called the Dilaton Equation (DE), the String Equation (SE)
and the Topological Recursion Relations (TRR), respectively.  The
Dilaton Equation implies that $\cL$ is a cone:

\begin{pro}
  \label{prop:cone}
  For every element $f\in \cL(R)$ and every $x\in R^{\rm nilp}$, we have
  $(1+x)f\in \cL(R)$.
\end{pro} 
\begin{proof} 
  Note that $f\in \cL(R)$ can be regarded as an element of
  $\cL(R[\![t]\!])$ by the natural inclusion $R \subset R[\![t]\!]$.
  It suffices to prove that $R[\![t]\!]$-valued point $(1+t) f \in
  (\cH,-z)(R[\![t]\!])$ belongs to $\cL(R[\![t]\!])$.  The conclusion
  then follows by applying the functor $\cL$ to the continuous
  $R$-homomorphism: $R[\![t]\!] \to R,\ t\mapsto x$.  We write $f = -z
  + h = -z + \bt(z) + \bp(z)$.  Then $(1+t)f = -z + h_t = -z +
  \bt_t(z) + \bp_t(z)$ with $h_t=-tz+(1+t)h$, $\bt_t(z) = -tz + (1+t)
  \bt(z)$ and $\bp_t(z) = (1+t)\bp(z)$.  Because $f\in \cL(R)$,
  $E_{j\beta}(h_t) = O(t)$.  Assume by induction on $n$ that
  $E_{j\beta}(h_t) = O(t^n)$.  We have
\begin{align*} 
&(1+t)\frac{d}{dt} E_{j\beta}(h_t) 
= (1+t) \left( 
p_{j\beta} -\sum_{k\ge 0} 
t_k^\alpha \corrr{\phi_\alpha \psi^k, \phi_\beta \psi^j}_{\bt_t} + 
\corrr{\psi,\phi_\beta \psi^j}_{\bt_t} \right) \\   
&\equiv \corrr{\phi_\beta\psi^j}_{\bt_t} 
- \sum_{k\ge 0} (-t \delta_k^1 \delta^\alpha_1 + (1+t)t_k^\alpha)
\corrr{\phi_\alpha \psi^k, \phi_\beta\psi^j}_{\bt_t} 
+ \corrr{\psi, \phi_\beta \psi^j}_{\bt_t} \mod t^n. 
\end{align*} 
In the second line, we used the induction hypothesis 
$E_{j\beta}(h_t) = O(t^n)$. 
But the second line is zero by the Dilaton Equation. 
Hence $E_{j\beta}(h_t) = O(t^{n+1})$. 
\end{proof}

The String Equation and the Topological Recursion Relations together
imply that the tangent space to $\cL$ at an $R$-valued point has the
structure of an $R\{z\}$-module:

\begin{pro} 
  \label{prop:closedunderz}
  The tangent space $T_f\cL(R)$ at $f\in \cL(R)$ is an
  $R\{z\}$-submodule of $T_f(\cH,-z)(R) = H^\bullet_{\rm
    orb}(\cX)\otimes R\{z,z^{-1}\}$.
\end{pro} 
\begin{proof} 
  Assume that we know $z T_f\cL(R) \subset T_f \cL(R)$. Then for every
  sequence $\{a_n\}_{n\ge 0}$ in $R$ with $\lim_{n\to \infty} a_n =0$
  and $h\in T_f\cL(R)$, we have
  \[
  a_0 h + a_1 z h + \cdots + a_m z^m h  \in T_f \cL(R).
  \] 
  This element converges to $(\sum_{n=0}^\infty a_n z^n)h$ as $m\to
  \infty$.  Since $T_f\cL(R)$ is a closed subspace,
  $(\sum_{n=1}^\infty a_n z^n)h \in T_f\cL(R)$.  Hence $T_f\cL(R)$ is
  an $R\{z\}$-submodule.

  Now it suffices to show that $z T_f\cL(R) \subset T_f \cL(R)$.  Take
  a tangent vector $h = \dot{\bt}(z) + \dot{\bp}(z)$ at $f$.  Then we
  have
  \begin{align*}
    z h = \sum_{k\ge 1} \dot{t}_{k-1}^\alpha \phi_\alpha z^k 
    - \dot{p}_{0\beta} \phi^\beta
    + \sum_{l\ge 0} (-\dot{p}_{l+1,\beta})
    \frac{\phi^\beta}{(-z)^{l+1}}.  
  \end{align*} 
  Therefore $zh \in T_f\cL(R)$ is equivalent to the equality:
  \[
  -\dot{p}_{l+1,\beta} = - \dot{p}_{0\nu}
  \corrr{\phi^\nu, \phi_\beta \psi^l}_{\bt} 
  + \sum_{k\ge 1} \dot{t}_{k-1}^\alpha 
  \corrr{\phi_\alpha \psi^k, \phi_\beta \psi^l }_{\bt}.  
  \]
  Substituting for $\dot{p}_{l\beta}$ using \eqref{eq:tangentvector},
  one sees that it suffices to show that 
  \begin{equation}
    \label{eq:k+1_l+1}
    \corrr{ \phi_\alpha \psi^{k+1}, \phi_\beta \psi^l}_{\bt} 
    + \corrr{\phi_\alpha \psi^k, \phi_\beta \psi^{l+1}}_{\bt}
    - \corrr{\phi_\alpha \psi^k, \phi_\nu}_{\bt}
    \corrr{\phi^\nu, \phi_\beta\psi^l }_{\bt} = 0.  
  \end{equation} 
  At $\bt=0$, we have (writing $\corrr{\cdots\vphantom{\psi^k}}_0$ for
  $\corrr{\cdots\vphantom{\psi^k}}_{\bt} \big \vert_{\bt = 0}$)
  \begin{align*}
    \corrr{\phi_\alpha \psi^{k+1}, \phi_\beta \psi^l}_0 +
    \corrr{\phi_\alpha \psi^k, \phi_\beta \psi^{l+1} }_0 &=
    \corrr{\phi_\alpha \psi^{k+1}, \phi_\beta\psi^{l+1}, \fun}_0 &&
    \text{by (SE)}  \\
    &= \corrr{\phi_\alpha \psi^k, \phi_\nu }_0 
    \corrr{\phi^\nu, \phi_\beta \psi^{l+1}, \fun}_0 &&
    \text{by (TRR)} \\
    & = \corrr{\phi_\alpha \psi^k, \phi_\nu}_0
    \corrr{\phi^\nu,\phi_\beta\psi^l}_0 && \text{by (SE)}.
\end{align*} 
On the other hand, differentiating in $t_j^\gamma$, we have
\begin{align*} 
  \partial_{j,\gamma} 
  &\left(\corrr{\phi_\alpha\psi^{k+1}, \phi_\beta \psi^l}_{\bt} 
    + \corrr{\phi_\alpha \psi^k, \phi_\beta \psi^{l+1} }_{\bt}\right) \\
  &= \corrr{\phi_\gamma \psi^j, \phi_\alpha\psi^{k+1},
    \phi_\beta\psi^l}_{\bt} 
  + \corrr{\phi_\gamma \psi^j, \phi_\alpha\psi^k,
    \phi_\beta \psi^{l+1}}_{\bt} \\ 
  &= \corrr{\phi_\alpha \psi^k, \phi_\nu}_{\bt} 
  \corrr{\phi^\nu, \phi_\gamma \psi^j, \phi_\beta \psi^l}_{\bt} 
  + \corrr{\phi_\beta \psi^l, \phi_\nu}_{\bt} 
  \corrr{\phi^\nu, \phi_\gamma\psi^j, \phi_\alpha\psi^k}_{\bt} 
  \\
  &= \partial_{j,\gamma} \left(
    \corrr{\phi_\alpha \psi^k, \phi_\nu}_{\bt}
    \corrr{\phi^\nu, \phi_\beta\psi^l }_{\bt}\right),  
\end{align*} 
where we used (TRR) in the third line.  Therefore we have
(\ref{eq:k+1_l+1}).
\end{proof} 

\begin{rem*}
Let $R_1,R_2$ be complete $\Lambda$-algebras. 
For a continuous $\Lambda$-algebra homomorphism 
$\varphi\colon R_1\to R_2$, the induced homomorphism 
$\varphi_* \colon T_{f}\cL(R_1) \to T_{\varphi(f)}\cL(R_2)$ 
becomes a continuous $R_1\{z\}$-module 
homomorphism. 
\end{rem*} 

Define elements $\tau^\alpha(\bt) \in S$, $1 \leq \alpha\leq N$, by
\[
\tau^{\alpha}(\bt) := 
\corrr{\fun, \phi^\alpha \vphantom{\psi^k}}_{\bt}.
\]
The String Equation implies that
\[
\tau^\alpha(\bt) = t_0^\alpha +  
\sum_{k\ge 0} t_{k+1}^\gamma 
\corrr{\phi_\gamma \psi^k, \phi^\alpha}_{\bt}
= t_0^\alpha + \text{higher order terms}.   
\]
Finally we establish the most remarkable property of $\cL$: that
tangent spaces to $\cL$ are parametrized by finitely many parameters
$\tau^1(\bt),\dots,\tau^N(\bt)$, and are generated by the derivatives
of the $J$-function as $R\{z\}$-modules.  This leads us to the
$\cD$-module property of tangent spaces (Corollary \ref{cor:Dmod}).

\begin{pro}
  \label{prop:tangentspace_derivativesofJ} 
  The tangent space $T_f\cL(R)$ at $f = -z +\bt(z) + \bp(z)\in \cL(R)$
  is freely generated by the derivatives of the $J$-function
  \[
  \partial_\alpha J(\tau,-z)|_{\tau=\tau(\bt)}, \quad \alpha=1,\dots,N. 
  \]
  as an $R\{z\}$-module, where $J(\tau,z)=J^\untw(\tau,z)$ is the
  untwisted $J$-function.
\end{pro} 

This is an immediate consequence of the following two lemmas.
\begin{lem} \label{thm:appendixlemma1} For elements $r^1,\dots,r^N\in
  R^{\rm nilp}$, the $J$-function $J(\tau,-z)$ with $\tau=\sum_\alpha
  r^\alpha\phi_\alpha$ gives an $R$-valued point on $\cL$.  The
  tangent space $T_{J(\tau,-z)}\cL(R)$ is freely generated by the
  derivatives $(\partial_\alpha J)(\tau,-z)$ as an $R\{z\}$-module.
\end{lem}

\begin{lem} \label{thm:appendixlemma2} The tangent space $T_f\cL(R)$
  at an $R$-valued point $f = -z + \bt(z) + \bp(z) \in \cL(R)$ is the
  same as the tangent space $T_{J(\tau(\bt),-z)} \cL(R)$ at
  $J(\tau(\bt),-z) \in \cL(R)$ as a subspace of $H_{\rm
    orb}(\cX)\otimes R\{z,z^{-1}\}$.  Here the ring homomorphism
  $\Lambda_\bs[\![\tau]\!] \to R$ sending $\tau^\alpha$ to
  $\tau^\alpha(\bt)\in R$ gives a point $J(\tau(\bt),-z) \in \cL(R)$.
\end{lem}

\begin{proof}[Proof of Lemma~\ref{thm:appendixlemma1}]
  As we saw earlier, the $J$-function is a
  $\Lambda[\![\tau]\!]$-valued point on $\cL$.  Thus its derivatives
  $\partial_\alpha J(\tau,-z)$ belong to the tangent space
  $T_{J(\tau,-z)}\cL(\Lambda_\bs[\![\tau]\!])$, by Lemma
  \ref{lem:diff_arc}.  Via the homomorphism $\Lambda[\![\tau]\!] \to
  R$ sending $\tau^\alpha$ to $r^\alpha \in R$, we obtain an
  $R$-valued point $J(\tau,-z) \in \cL(R)$ and tangent vectors
  $(\partial_\alpha J)(\tau,-z) \in T_{J(\tau,-z)}\cL(R)$.

  Set $\tau = \sum_{\alpha} r^\alpha \phi_\alpha$ and write $[f]_+$
  for the non-negative part of the $z$-series $f$.  From the
  description (\ref{eq:tangentvector}) of tangent vectors, there
  exists a one-to-one correspondence between tangent vectors
  $\dot{\bt}(z) + \dot{\bp}(z)$ in $T_{J(\tau,-z)} \cL(R)$ and tuples
  $\{\dot{t}_k^\alpha\in R\}_{k,\alpha}$ satisfying $\lim_{k\to
    \infty} \dot{t}_k^\alpha =0$.  It therefore suffices to show that
  for any given $\{\dot{t}_k^\alpha\in R\}_{k,\alpha}$ satisfying
  $\lim_{k\to\infty} \dot{t}_k^\alpha =0$, there exist unique elements
  $c^\alpha \in R\{z\}$ such that
  \begin{equation}
    \label{eq:derivatives_J_generate}
    \left[\sum_{\alpha} c^\alpha 
      \partial_\alpha J \right]_+ =  
    \sum_{k\ge 0} \dot{t}_k^\alpha \phi_\alpha z^k.  
  \end{equation} 
  First we show the existence of $c^\alpha$.  Assume that the topology
  on $R$ is defined by a descending chain of ideals $\{I_M\}_{M\ge
    0}$.  We will prove the following claim by induction on $n$: 

  \begin{cla*}
    There exist $c(n)^\alpha \in R\{z\}$ such that
    \[
    \left[
      \sum_{\alpha} c(n)^\alpha \partial_\alpha J
    \right ]_+ = \sum_{k=0}^n 
    \dot{t}_k^\alpha \phi_\alpha z^k.  
    \]
    Moreover, if $\dot{t}_k^\beta\in I_M$ for all 
    $0\le k\le n$ and $\beta$ and for some $M$,  
    then $c(n)^\alpha \in I_M\{z\}$.
  \end{cla*}
  
  The case $n=0$ is clear from the expansion
  \[
  \partial_\alpha J(\tau,-z) = \phi_\alpha + 
  \sum_{j\ge 0} \corrr{\phi_\alpha, \phi_\beta\psi^j}_\tau 
  \frac{\phi^\beta}{(-z)^{j+1}}.  
  \]
  Assume that the claim holds for some $n\ge 0$.  One then has
  \begin{align*}
    \left[
      \sum_\alpha 
      (\dot{t}^\alpha_{n+1}z^{n+1} + c(n)^\alpha) 
      \partial_\alpha J 
    \right]_+ 
    = \sum_{k=0}^{n+1} \dot{t}_k^\alpha \phi_\alpha z^k 
    - \sum_{j=0}^n \dot{t}_{n+1}^\alpha 
    \corrr{\phi_\alpha,\phi_\beta\psi^j}_{\tau} 
    (-1)^{j} \phi^\beta z^{n-j}.  
  \end{align*}
  By the induction hypothesis, there exist $\xi(n)^\alpha\in R\{z\}$
  such that $[\sum_\alpha \xi(n)^\alpha \partial_\alpha J]_+ =
  \sum_{j=0}^n \dot{t}_{n+1}^\alpha
  \corrr{\tau_{0\alpha}\tau_{j\beta}}_{\tau} (-1)^{j} \phi^\beta
  z^{n-j}$.  Also we have $\xi(n)^\alpha \in I_M\{z\}$ if
  $\dot{t}_{n+1}^\beta\in I_M$ for all $\beta$.  Therefore, we can
  take $c(n+1)^\alpha$ to be $c(n)^\alpha + \dot{t}_{n+1}^\alpha
  z^{n+1} + \xi(n)^\alpha$.  This completes the induction step, and
  the claim follows.

  The above argument shows that $c(n+1)^\alpha-c(n)^\alpha \in
  I_M\{z\}$ if $\dot{t}_{n+1}^\beta\in I_M$ for all $\beta$.
  Therefore $c(n)^\alpha$ converges to some element $c^\alpha\in
  R\{z\}$ and (\ref{eq:derivatives_J_generate}) holds.  For the
  uniqueness of $c^\alpha$, it suffices to show that if $[\sum_\alpha
  c^\alpha \partial_\alpha J]_+=0$ then $c^\alpha = 0$. Suppose that
  $[\sum_\alpha c^\alpha \partial_\alpha J]_+=0$ and that $c^\beta
  \neq 0$ for some $\beta$. Since $R$ is Hausdorff, there exists an
  $M$ such that $c^\beta \notin I_M\{z\}$.  The equation $[\sum_\alpha
  c^\alpha \partial_\alpha J]_+=0$ holds in the ring
  $R\{z,z^{-1}\}/I_M\{z,z^{-1}\} = (R/I_M)[z,z^{-1}]$ and $c^\beta
  \neq 0$ in $(R/I_M)[z]$.  Comparing the highest order terms in $z$
  leads us to a contradiction.
\end{proof} 

\begin{proof}[Proof of Lemma~\ref{thm:appendixlemma2}]
  Recall again that a tangent vector at $f$ is given by a set
  $\{\dot{t}_k^\alpha, \dot{p}_{k\alpha}\}_{k,\alpha}$ in $R$
  satisfying $\lim_{k\to \infty} \dot{t}_k^\alpha = \lim_{k\to \infty}
  \dot{p}_{k\alpha} =0$ and equation \eqref{eq:tangentvector}.  On the
  other hand, the Topological Recursion Relations imply that:
  \[
  \corrr{\phi_\alpha \psi^k, \phi_\beta \psi^l}_{\bt} 
  = \corrr{\phi_\alpha \psi^k, \phi_\beta \psi^l}_{\tau(\bt)}.  
  \] 
  This is due to Dijkgraaf--Witten \cite{Dijkgraaf-Witten} (see also
  \cite[Equation~2]{Givental:symplectic},
  \cite[Proposition~4.6]{Getzler:jetspace}). The Lemma follows.
\end{proof}

\begin{cor}
  \label{cor:Dmod} 
  Let $I(t)$ be an $R[\![t]\!]$-valued point on $\cL$ and $\xi(t)$ be a
  tangent vector at $I(t)$. Then $z {d\xi \over dt}(t)$ is again a
  tangent vector at $I(t)$.
\end{cor} 
\begin{proof}
  Set $I(t) = -z + \bt_t(z) + \bp_t(z)$. 
  By Proposition \ref{prop:tangentspace_derivativesofJ}, 
  we can write $\xi(t)$ in the form 
  \begin{align*}
    \xi(t) = \sum_{\alpha} c^\alpha(t,z) [\partial_\alpha
    J(\tau,-z)]_{\tau=\tau(\bt_t)}, &&
    \text{for some $c^\alpha(t,z) \in R[\![t]\!]\{z\}$.}
  \end{align*}
  The Corollary follows from this and the differential equations
  (\ref{eq:qdes}). 
\end{proof}

\end{document}